\documentclass[a4paper,12pt]{article}
\usepackage{latexsym,amsthm,hhline}
\newtheorem{define}{Definition}
\newtheorem{lemma}{Lemma}

\newtheorem{theorem}{Theorem}

\def \H {{I\!\!H}}
\def \E {{I\!\!E}}
\def \Z {{Z\!\!\!Z}}
\makeatletter
\def\mmcaption #1{\refstepcounter{table}
\begingroup
  \@parboxrestore
  \normalsize
  \@makecaption{\fnum@table}{\ignorespaces #1}\par
\endgroup}

\newenvironment{mmtable}
              {\medskip\vbox\bgroup}
              {\egroup\medskip}
\makeatother

\begin{document}

\begin{center}
{\LARGE\bf
Coxeter decompositions of hyperbolic simplices.   }

\medskip
{\Large A.~Felikson}
\end{center}

\setcounter{section}{0}
\section{Introduction}
\begin{define}
A convex polyhedron in a space of constant curvature
is called {\bf a Coxeter polyhedron}
if all  dihedral angles of this polyhedron are the integer parts of
$\pi$.

\end{define}

\begin{define}
\label{def1}
{\bf A  Coxeter decomposition} of a convex polyhedron $P$
is a decomposition of $P$ into finitely many tiles
such that each tile is a
Coxeter polyhedron and
any two tiles having a common facet are symmetric
 with respect to this facet.
\end{define}

Coxeter decompositions of hyperbolic triangles were studied
in~\cite{Pink}, \cite{K}, \cite{Mat}, \cite{M} and \cite{Deza}.
Coxeter decompositions of hyperbolic tetrahedra are listed
in~\cite{hyp3}.
In this paper, we classify Coxeter decompositions of simplices
in hyperbolic spaces $\H^n$, where $n\ge 4$.
The paper completes the classification of the Coxeter decompositions
of hyperbolic simplices.

The author is grateful to E.~B.~Vinberg for the numerous very helpful
remarks.

\subsection*{Basic definitions}

The tiles in  Definition~\ref{def1} are called
{\bf fundamental polyhedra}.
Clearly, any two fundamental polyhedra are congruent to each other.
A hyperplane $\alpha$ containing a facet of a fundamental polyhedron
is called a {\bf mirror} if $\alpha$ contains no facet of $P$.

\begin{define}
Given a Coxeter decomposition of a polyhedron $P$,
a {\bf dihedral angle} of $P$ formed up by facets
$\alpha$ and $\beta$ is called {\bf fundamental}
if no mirror contains $\alpha \cap \beta$.
\end{define}

From now on we consider only the polyhedra
bounded by the mirrors of some Coxeter decomposition.

\medskip

\noindent
{\bf Notation.}

\noindent
$P$ is a simplex equipped with a Coxeter decomposition; \\
$F$ is a fundamental polyhedron
considered up to an isometry of $\H^n$; \\
$\Sigma(T)$ is a Coxeter diagram of a Coxeter simplex  $T$; \\
$N$ is a number of the fundamental polyhedra inside $P$.\\

\noindent
A decomposition is called {\bf non-trivial} if $N>1$.

\medskip

\noindent

We use the standard notation for the Euclidean and spherical
Coxeter simplices: $A_n$,
$\widetilde A_n$,
$B_n$,
$\widetilde B_n$,
$\widetilde C_n$,
$D_n$,
$\widetilde D_n$,...
The notation for the hyperbolic Coxeter simplices is introduced in
Table~2.

\section{Properties of Coxeter decompositions \\ of simplices}

\subsection{Fundamental polyhedron}

\begin{theorem}
If $P$ is a hyperbolic simplex then $F$ is a simplex too.

\end{theorem}

\begin{proof}
Since the decomposition contains finitely many tiles,
there is a {\bf minimal } simplex inside $P$
containing no proper simplices bounded by mirrors and faces of $P$.
We are aimed to prove that the minimal simplex is a fundamental polyhedron
of the decomposition.
Without loss of generality, we assume that $P$ is minimal itself.

Suppose that $P$ has a non-fundamental dihedral angle.
Let $\Pi$ be a mirror decomposing the dihedral angle of $P$.
Then $\Pi$ decomposes $P$ into two smaller simplices.
The contradiction to the minimal property of $P$ shows
that any dihedral angle of $P$ is fundamental.

Suppose that $F$ is not a simplex.
Consider the set $S$ of polyhedra which could be cut from $P$
by a single mirror. The number of mirrors is finite.
Thus, $S$ contains a {\bf minimal}  polyhedron $M$
such that no polyhedron contained in $M$ belongs to $S$.
Let $\Pi$ be a mirror which cuts $M$ from $P$.

Consider the dihedral angles of $M$.
Any dihedral angle $\alpha$ of $M$ formed up by $\Pi$
and a facet of $P$ is fundamental (otherwise
$M$ contains a polyhedron which belongs to $S$
in contradiction to the minimal property of $M$).
The rest dihedral angles of $M$ are the dihedral angles of $P$.
Thus, all dihedral angles of $M$ are fundamental,
and $M$ is a Coxeter polyhedron.

Let $A_0$,...,$A_k$ be the vertices of $P$ contained in
$M$, let $A_{k+1}$,...,$A_n$ be the rest vertices of $P$.
If $k=0$ then $M$ is a simplex in contradiction to the minimal property of
$M$. Therefore, $k\ne 0$. Consider the lines $A_0A_n$ and $A_1A_n$.
These lines contain the non-adjacent edges of $M$.
But these lines have a common point $A_n$.
This is impossible, since $M$ is an acute angled polyhedron
(see~\cite{Andr}).
The contradiction shows that $F$ is a simplex.

\end{proof}

\begin{lemma}
\label{fin}
Let $F$ be a hyperbolic Coxeter simplex.
The number of simplices admitting a Coxeter decomposition with
the fundamental polyhedron $F$ is finite.

\end{lemma}

\begin{proof}
The dihedral angles completely determine a non-Euclidean simplex.
The dihedral angles of a simplex $P$ admitting a Coxeter decomposition
are the multiples of the dihedral angles of the fundamental polyhedron.
Thus, for any fundamental simplex $F$ there exist finitely many
ways to prescribe the dihedral angles of the simplex $P$.

\end{proof}

\begin{lemma}[Volume property]
$\frac{Vol(P)}{Vol(F)}\in \mathbf Z$, where $Vol(T)$ is a volume
of a simplex $T$.

\end{lemma}

The lemma is evident.
The volumes of the hyperbolic Coxeter simplices are computed in~\cite{Ruth}.
In Table~2, we list the volumes of the Coxeter simplices in $\H^n$,
where $n \ge 4$.

\subsection{Three types of Coxeter decompositions of simplices}
\label{ind}

Consider a Coxeter decomposition of a simplex $P$.
Suppose that $P$ has a non-fundamental dihedral angle.
Suppose in addition, that any non-fundamental simplex inside $P$
has a non-fundamental dihedral angle.  Let $\Pi$ be a
mirror decomposing the dihedral angle of $P$.  Then $\Pi$ decomposes $P$
into two smaller simplices $P_1$ and $P_2$.  Given the decompositions of
$P_1$ and $P_2$, it is possible to find the decomposition of $P$.  By the
assumption each of $P_1$ and $P_2$ is either fundamental or has a
non-fundamental dihedral angle.  In the latter case, the small simplex is
decomposed into two smaller ones.  Thus, we decompose the simplex into
smaller and smaller simplices.  The process stops only when $P$
is decomposed into fundamental simplices.  Since $P$ contains finitely
many fundamental simplices, the process stops after finitely many steps.

Invert the process. We obtain the following inductive algorithm
for constructing of the decomposition of the big simplex from the
 decompositions of the smaller simplices:

\begin{itemize}
\item[]
\noindent
$\bullet$
Step 0.
Let $F$ be a fundamental simplex. Set $P_0=F$.\\
\noindent
$\bullet$
Step 1.
Take two copies of $P_0$.
Glue these copies together to construct a simplex admitting
a Coxeter decomposition with $N=2$.
Use all the ways to glue these copies together and
 obtain all the simplices $P_1,...,P_k$ admitting the decomposition
with this fundamental polyhedron and $N=2$.\\
\noindent
$\bullet$
Step 2.
Find all simplices
$P_{k+1},...,P_l$ which consist either of $P_1$ and $P_0$
or of two copies of $P_1$
(in the former case we take one copy of each simplex).
\\
\phantom{}  . . . . . . . . . . . .
. . . . . . . . . . . . . . . . . . . . . . . . . . . . . . .
\\
\noindent
$\bullet$
Step m.
Suppose that after the step $(m-1)$
we have the simplices $P_0,P_1,...,P_n$.
At the step $m$ compose all the simplices
$P_{n+1},...,P_s$
from
one copy of
$P_{m-1}$ and one copy of $P_i$ for each $i\le m-1$
(if the simplex obtained is already in the list  $P_0,P_1,...,P_n$,
one should not add this simplex to the list).

\end{itemize}

Suppose that $P$ admits a decomposition with the fundamental simplex $F$.
Suppose that any non-fundamental simplex inside $P$
has a non-fundamental dihedral angle.
Then the decomposition of $P$ will be constructed after finitely many
steps of the inductive algorithm.
In general, a simplex constructed by the algorithm contains some
non-fundamental simplices with fundamental dihedral angles.
But any simplex obtained by the algorithm has a non-fundamental dihedral
angle.

\begin{define}
Let $\Theta(P)$ be a Coxeter decomposition of a simplex $P$.

$\Theta(P)$ is called a decomposition of the
{\bf first type} if the decomposition can be obtained by the inductive
algorithm.

$\Theta(P)$ is called a decomposition of the
{\bf second type} if all the dihedral angles of $P$ are fundamental.

$\Theta(P)$ is called a decomposition of the
{\bf third type} if $\Theta(P)$ is neither of the first type nor of the
second.

\end{define}

\medskip

To find the decompositions of the first type, we use the inductive
algorithm. To classify  the decompositions of the second type is more
difficult, since there is no algorithm for the classification.
It turns out that any decompositions of the third type is a superposition
of the decompositions of the first and the second types.

\section{Classification of the decompositions}

\subsection{Decompositions of the first type}
\label{type1}

Consider the inductive algorithm for the fundamental simplex $F$.
The number of simplices admitting a decomposition with the fundamental
simplex $F$ is finite.
Each decomposition can be obtained by finitely many ways.
Thus, the inductive algorithm terminates after finitely
many steps.
The result of the algorithm is the complete set of the
decompositions of the first type with the fundamental simplex $F$.

There are finitely many hyperbolic Coxeter simplices in the dimensions
greater than 2. Thus, we can use the inductive algorithm for one Coxeter simplex
after another. We can use a computer for the computations.

\begin{define}
Let $\Theta (F,P)$ be a Coxeter decomposition of the simplex $P$
with fundamental simplex $F$.
Let $\Theta _1(T,P)$ be a Coxeter decomposition of the simplex $P$
with fundamental simplex $T$.
Suppose that each mirror of  $\Theta _1(T,P)$
is a mirror of $\Theta (F,P)$.
Then $\Theta (F,P)$  is called a {\bf superposition}
of the decompositions $\Theta _1(T,P)$ and $\Theta (F,T)$,
where $\Theta (F,T)$ is a restriction of $\Theta (F,P)$ to the simplex $T$.
A decomposition is called {\bf simple}
if it is not a superposition of the non-trivial decompositions.

\end{define}

Clearly, it is sufficient to list the simple decompositions.
Tables 3 and 4 contain the list of the simple decompositions of
simplices of the first type.

\subsection{Decompositions of the second type} \label{type2}

Let $P$ be a simplex admitting a Coxeter decomposition of the second type.
Clearly, $P$ is a Coxeter simplex.

\begin{lemma}[Subdiagram property]
\label{subd_3}
Let $P$ be a simplex admitting a Coxeter decomposition of the second type
with fundamental simplex $F$.
Let $\Sigma(P)$ be a Coxeter diagram of $P$ and $\Sigma(F)$ be a Coxeter
diagram of $F$.  Let $v$ be a node of $\Sigma(P)$.
There exists a node $w$ of $\Sigma(F)$ such that
either
$\Sigma(P)\setminus v = \Sigma(F) \setminus w$
or the simplex  $p$ determined by the subdiagram $\Sigma(P) \setminus v$
admits a Coxeter decomposition of the second type with fundamental simplex
$f$ determined by $\Sigma(F) \setminus w$.

\end{lemma}

\begin{proof}
The node $v$ of $\Sigma (P)$ corresponds to a facet $\alpha$ of $P$.
Let $A$ be a vertex of $P$ opposite to the facet $\alpha$.

Suppose that the fundamental simplex $F_A$ containing $A$ is unique.
Then all but one facets of $F_A$ are the facets of $P$.
Let $w$ be a node of $\Sigma (F)$ correspondent to the rest facet of
$F_A$. Then  $\Sigma(P) \setminus v =
\Sigma(F) \setminus w$.

Suppose that $A$ belongs to several fundamental simplices.
Consider a small sphere $s$ centered in $A$ (if the vertex $A$ is ideal,
consider the horosphere $s$ centered in $A$).
Let $p=P\bigcap s$.
The Coxeter decomposition of $P$ restricted to  $p$
is a Coxeter decomposition of a spherical (Euclidean)
simplex $p$ with some fundamental simplex $f$.
Clearly, a  Coxeter diagram
of $p$ is $\Sigma(P) \setminus v$, and a Coxeter diagram
of $f$ is $\Sigma(F) \setminus w$ for some face $w$ of $F$.
Since the dihedral angles of $P$
are fundamental, the dihedral angles of $p$ are fundamental too.

\end{proof}

To use the Subdiagram property we need the classification of the
Coxeter decompositions of spherical simplices of the second type
(see~\cite{sph_sh}).

\begin{define}
A Coxeter decomposition of a spherical polyhedron is called
{\bf indecomposable} if the Coxeter diagram of the fundamental simplex is
connected.

\end{define}

In~\cite{sph_sh}, it is proved that any Coxeter decomposition of
a spherical simplex is a direct product of the indecomposable decompositions.
Moreover, a decomposition of the second type is a product of the
indecomposable  decompositions of the second type.

\begin{theorem}\label{sphere2}
An indecomposable decomposition of a spherical simplex of the second type
is uniquely determined by the pair $(F,P)$.
The possibilities for the pairs $(F,P)$ are listed in Table~\ref{sph_tab}.

\end{theorem}

\begin{mmtable}
\begin{center}
\mmcaption{Indecomposable Coxeter decompositions}
of spherical simplices of the second type.
\label{sph_tab}
\vspace{13pt}
{\begin{tabular}{|c|c|}
\hline
$\vphantom{\int^A_A}$
$F$
$\vphantom{\int^A_A}$
&$P$\\
\hhline{|=|=|}
$\vphantom{\int^A_A}$
$H_3$
$\vphantom{\int^A_A}$
&$3A_1$\\
\hline
$\vphantom{\int^A_A}$
$F_4$
$\vphantom{\int^A_A}$
&$2A_2$\\
\hline
$\vphantom{\int^A_A}$
$H_4$
$\vphantom{\int^A_A}$
&
\begin{tabular}{c}
$A_4$, $2G_2^{(5)}$, $2A_2$,
$H_3+A_1$, $D_4$, $4A_1$
\end{tabular}\\
\hline
$\vphantom{\int^A_A}$
$D_n$&
$\vphantom{\int^A_A}$
\begin{tabular}{l}
$D_{m_1}+...+D_{m_r}$ ($m_1+...+m_r=n$),\\
$m_1\ge m_2\ge ... \ge m_r>1$, where $D_2=2A_1$, $D_3=A_3$
\end{tabular}\\
\hline
$\vphantom{\int^A_A}$
$E_6$
$\vphantom{\int^A_A}$
&$A_5+A_1$, $3A_2$\\
\hline
$\vphantom{\int^A_A}$
$E_7$
$\vphantom{\int^A_A}$
&$D_6+A_1$, $A_5+A_2$,
$2A_3+A_1$, $A_7$, $D_4+3A_1$,
$7A_1$\\
\hline
$\vphantom{\int^A_A}$
$E_8$
$\vphantom{\int^A_A}$
&
\begin{tabular}{c}
$A_8$,
$A_7+A_1$,
$A_5+A_2+A_1$,
$2A_4$,
$4A_2$, $A_6+A_2$,
$E_7+A_1$,\\
$D_8$,
$D_6+2A_1$,
$D_5+A_3$, $2D_4$,
$D_4+4A_1$, $2A_3+2A_1$, $8A_1$

\end{tabular}
\\
\hline

\end{tabular}

 }

\end{center}
\end{mmtable}

If $F\ne H_3,H_4,F_4,B_n$,
the classification follows immediately from~\cite{Dyn}.
The rest cases of Theorem~\ref{sphere2} are proved in~\cite{sph_sh}.

\begin{lemma} \label{versh}
Let $P$ be an $n$-dimensional simplex and $\Theta(P)$
be a non-trivial Coxeter decomposition such that
for any vertex $v$ of $P$ a fundamental simplex $F_v$ containing
$v$ is unique.
Then $v$ is the only vertex of $P$ contained in $F_v$.

\end{lemma}

\begin{proof}
Let $\alpha _v$ be a facet of $F_v$ opposite to $v$.
The rest facets of $F_v$ are the facets of $P$.
Suppose that  $F_v$ contains a vertex $w\ne v$ of $P$.
Then $\alpha _v$ contains $w$.
Since the decomposition is non-trivial, $\alpha _v$ is not a facet of $P$.
Thus, $w$ belongs to at least two fundamental simplices,
that is impossible by the hypothesis of the lemma.

\end{proof}

\begin{lemma} \label{8}
Let $P$ be an $n$-dimensional simplex and $\Theta(P)$
be a Coxeter decomposition such that
for any vertex $v$ of $P$ a fundamental simplex $F_v$ containing
$v$ is unique.
Then
 $N\ge 2(n+1)$.

\end{lemma}

\begin{proof}
Let $v$ be a vertex of $P$ and $F_v$ be the fundamental simplex
containing $v$. Let $F^v$ be a fundamental simplex having a common facet
with $F_v$ (the fundamental simplex with this property is unique).
Let $F_w$ and $F^w$ be the analogous simplices for the vertex $w\ne v$.
To prove the lemma, it is sufficient to show that the simplices
$F_v$, $F^v$, $F_w$ and $F^w$ are distinct.

By the hypothesis of the lemma
$F_v\ne F^v$ and $F_w\ne F^w$.
Suppose that either $F_v =F_w$, or $F_v =F^w$,
or $F^v =F_w$, or $F^v =F^w$.
Then the simplices
$F_v\bigcap F^v$ and $F_w\bigcap F^w$
are the facets of the same $n$-dimensional simplex $S$.
By Lemma~\ref{versh},
$F_v\bigcap F^v$ has a vertex in the inner part of each edge of $P$
incident to $v$. The simplex
$F_w\bigcap F^w$ has a vertex in the inner part of each edge incident to
$w$. Thus, $S$ has at least $2n-1$ vertices.
Since $n>2$, this is impossible.

\end{proof}

\begin{lemma}  \label{2}
Let $\Theta (P)$ be a Coxeter decomposition containing exactly two
fundamental simplices. Then $P$ has a non-fundamental dihedral angle.

\end{lemma}

The lemma is evident.

\begin{lemma} \label{similar}
Let $P$ be a simplex in $\E^n$ admitting a non-trivial Coxeter
decomposition. Suppose that $F$ is similar to $P$.
Suppose that the Coxeter diagram
$\Sigma(F)$ differs from $\widetilde '_2$ and $\widetilde G_2^{(5)}$.
Then $N\ge 2^n$.

\end{lemma}

\begin{proof}
Let $\Sigma (P)$ be the affine Dynkin diagram
and $\bar v$ be the node of $\Sigma (P)$ correspondent to the lowest
root.  Let $v$ be the vertex of $P$ opposite to the facet that
is represented by $\bar v$. Let $\widetilde W$ be the affine Weyl group.
The stabilizer of $v$ in $\widetilde W$ is the set of all linear parts of
the isometries contained in $\widetilde W$.  Thus, the fundamental simplex
$F$ containing $v$ is unique.  Moreover, if $\Sigma(F)$ differs from
$\widetilde '_2$ and $\widetilde G_2^{(5)}$ then $F$ is homothetic to
 $P$.

Suppose that $\Sigma(P)\ne \widetilde A_n$.
Let $\alpha $ be the facet of $F$ opposite to $v$.
Consider an edge $r_F$ of $F$ incident to $v$ and orthogonal to $\alpha$.
Let $r_P$ be an edge of $P$ containing $r_F$.
The reflection with respect to $\alpha$ preserves the line
containing $r_F$.
Hence, $r_P$ is at least to times longer than $r_F$,
and the lemma is proved.

Suppose that $\Sigma(P)=\widetilde A_n$.
To show that in this case the coefficient of the homothety is
at least 2, consider a set $R_{min}$ of the shortest edges of
$F_v$. The symmetry group of $F_v$ acts transitively on the
set of vertices of $F_v$.
Thus, $R_{min}$  contains at least one edge $r_v$ incident to $v$.
Let $t_v$ be the edge of $P$ containing $r_v$.
The fundamental simplex contains no edge shorter than $r_v$.
Hence, $t_v$ is at least to times longer than $r_v$,
and
$\frac{Vol P}{Vol F}\ge 2^n$.

\end{proof}

Let $(F,P)$ be a pair of simplices in $\H^n$, $n\ge 4$.
Suppose that $(F,P)$ satisfies  Volume Property,
 Subdiagram Property, Lemma~\ref{8} and Lemma~\ref{2}.
It turns out that in this case $(F,P)$ is one of the following pairs:\\
\phantom{qweq}
$(H^4_3,H^4_9),\ N=10$;\qquad
$(H^5_5,H^5_{12}),\ N=16$;\qquad
$(H^5_7,H^5_{11}),\ N=6$;\\
\phantom{qweq}
$(H^5_4,H^5_{11}),\ N=20$;\qquad
$(H^8_1,H^8_4),\ N=272$;\qquad
$(H^9_1,H^9_3),\ N=527$  \\
(see Table~2 for the notation of the hyperbolic Coxeter simplices).

Consider these five cases to figure out if there exist decompositions
corresponding to these pais $(F,P)$.
We call $(F,P)$,  the {\bf symbol} of the decomposition.
Sometimes we write the number $N$ near the pair $(F,P)$.

\subsubsection{There exists a unique decomposition with
symbol $(H^4_3,H^4_9)$ and $N=10$. }

The fundamental simplex  $H^4_3$ has a unique ideal vertex $v$.
The section of $H^4_3$  by the horosphere centered in $v$
is determined by the parabolic diagram $\widetilde A_3$
(we call $v$ the vertex of the type
$\widetilde A_3$).
The simplex $H^4_9$ has three ideal vertices of the type $\widetilde A_3$.
By Lemma~\ref{similar},
for any decomposition with symbol
$(\widetilde A_3,\widetilde A_3)$ we have
$N=1$ or $N\ge 8$.
Since 10 fundamental simplices have 10 ideal vertices,
one of the ideal vertices of $H^4_9$
(say, $v_1$) belongs to 8 fundamental simplices,
the rest two ideal vertices of $H^4_9$ belongs to one fundamental simplex
each.

The only decomposition with symbol $(\widetilde A_3,\widetilde A_3), N=8$
is a decomposition shown in
Fig.~1.
This is the decomposition of the section by horosphere centered in $v_1$.
Let $F_1$,...,$F_8$ be the fundamental simplices containing $v_1$,
let $f_1$,...,$f_8$ be the facets of $F_1$,...,$F_8$ opposite to $v_1$.
Let $S$ be the union of $f_1$,...,$f_8$.
The combinatorics of $S$ coincides with the decomposition shown
in Fig.~1.
But the facets  $f_1$,...,$f_8$ belong to different hyperplanes.
Indeed, $f_1$ is orthogonal to all but one facets of $F_1$.
The rest facet intersects $f_1$ with the dihedral angle $\frac{\pi}{3}$.
A 2-dimensional face of the fundamental simplex $H^4_3$
assigned with the dihedral angle $\frac{\pi}{n}$
will be called an {\bf $n$-face}.
Let $t$ be a 2-dimensional face of the fundamental simplex $\widetilde A_3$
in the decomposition with symbol
$(\widetilde A_3,\widetilde A_3)$.
We call $t$ an $n$-face if $t$ corresponds to the $n$-face of $H^4_3$.

It is not evident which face of the tetrahedron $AEFG$ is a 3-face
(see Fig.~1).

\begin{figure}
\label{e3}
\begin{center}
\caption{}
\vspace{10pt}
\setlength{\unitlength}{0.075in}
\begin{picture}(32.15,41.01)
\special{em:linewidth 0.014in}
\put(20.81,19.50){\special{em:moveto}}
\put(13.61,8.16){\special{em:lineto}}
\put(29.61,7.33){{\setbox0=\hbox{$M$}\lower\ht0\box0}}
\put(13.35,17.33){{\setbox0=\hbox{$L$}\lower\ht0\box0}}
\put(12.15,7.33){{\setbox0=\hbox{$K$}\lower\ht0\box0}}
\put(22.01,19.50){{\setbox0=\hbox{$G$}\lower\ht0\box0}}
\put(24.81,27.50){{\setbox0=\hbox{$F$}\lower\ht0\box0}}
\put(6.68,27.83){{\setbox0=\hbox{$E$}\lower\ht0\box0}}
\put(26.41,0.00){{\setbox0=\hbox{$D$}\lower\ht0\box0}}
\put(31.88,17.50){{\setbox0=\hbox{$C$}\lower\ht0\box0}}
\put(0.01,14.66){{\setbox0=\hbox{$B$}\lower\ht0\box0}}
\put(14.55,41.00){{\setbox0=\hbox{$A$}\lower\ht0\box0}}
\put(15.48,14.83){\special{em:moveto}}
\put(20.68,19.50){\special{em:lineto}}
\put(20.81,19.66){\special{em:moveto}}
\put(29.08,7.83){\special{em:lineto}}
\put(8.81,26.83){\special{em:moveto}}
\put(15.75,14.83){\special{em:lineto}}
\put(23.61,26.83){\special{em:lineto}}
\put(23.75,26.83){\special{em:moveto}}
\put(20.81,19.50){\special{em:lineto}}
\put(8.68,26.83){\special{em:lineto}}
\put(13.48,8.16){\special{em:moveto}}
\put(15.61,14.83){\special{em:lineto}}
\put(29.08,7.83){\special{em:lineto}}
\put(8.68,26.83){\special{em:moveto}}
\put(23.61,26.83){\special{em:lineto}}
\put(29.08,7.83){\special{em:lineto}}
\put(13.48,8.16){\special{em:lineto}}
\put(8.68,26.83){\special{em:lineto}}
\put(15.35,38.83){\special{em:moveto}}
\put(32.15,15.00){\special{em:lineto}}
\put(2.28,14.83){\special{em:moveto}}
\put(26.01,0.83){\special{em:lineto}}
\put(15.35,38.83){\special{em:moveto}}
\put(2.28,14.83){\special{em:lineto}}
\put(17.21,14.91){\special{em:lineto}}
\put(32.15,15.00){\special{em:lineto}}
\put(26.01,0.83){\special{em:lineto}}
\put(15.35,38.83){\special{em:lineto}}
\end{picture}
\end{center}
\vspace{15pt}
\end{figure}

\begin{itemize}
\item[]

\noindent
{\bf Proposition.}
Exactly two of 3-faces cut $ABCD$.
These facets are  $LMF$ and $LEK$.

\begin{proof}
The 3-faces cut the set $S$ into parts contained in different hyperplanes.
Only one of these hyperplanes (say, $\Pi$) contains a facet of
$H^4_9$.
Let $f_i\notin \Pi$. Then $f_i$ belongs to exactly two fundamental
simplices. Since the decomposition contains 10 fundamental simplices,
exactly two of the faces
$f_1$,...,$f_8$ are not in $\Pi$.
Thus, at most two 3-faces cut $S$.
Hence, at most two 3-faces correspond to the triangles in the inner part
of the tetrahedron $ABCD$ (see Fig.~1).
If a face is symmetric to another face with respect to some
mirror, then these faces are either 2-faces both or 3-faces both.
These conditions are satisfied only if
$KLM$, $LMF$, $LEF$, $LEK$, $KMD$ and $AEF$ are 3-faces,
and the rest triangles are 2-faces.
The 3-faces inside $ABCD$ are
$LMF$ and $LEK$.

\end{proof}
\end{itemize}

Consider a union of $F_1$,...,$F_8$.
Attach to the union of $F_1$,...,$F_8$ two additional simplices,
glue these simplices to the facets correspondent to the tetrahedra
$LMGC$ and $BKEL$.
We obtain a simplex
$H^4_9$ and the decomposition with symbol
$(H^4_3,H^4_9),N=10$.
Clearly, the decomposition with this symbol is unique.

\subsubsection{There is no decomposition with symbol
$(H^5_5,H^5_{12}),\ N=16$.
}
\label{type2.4}

The Coxeter diagram of $H^5_{12}$ has two parabolic subdiagrams
of the type
$\widetilde C_4$ and four parabolic subdiagrams of the type
$\widetilde F_4$.
The Coxeter diagram of $H^5_5$ has one subdiagram of the type
$\widetilde C_4$ and one subdiagram of the type $\widetilde F_4$.
A diagram $\widetilde F_4$ contains no subdiagram $2B_2$
and a diagram $\widetilde C_4$ contains no subdiagram $2A_2$.
Moreover, $2A_2$ and $2B_2$ admit no decomposition of the second type
with $N<32$.
Therefore, the vertex of $P$ of the type $\widetilde C_4$
is tiled by sixteen vertices of the type $\widetilde C_4$;
the vertex of the type $\widetilde F_4$ is tiled by
sixteen vertices of the type $\widetilde F_4$.
By Lemma~\ref{similar}, a non-trivial decomposition with symbol
($\widetilde F_4,\widetilde F_4$) contains at least 16 fundamental
simplices.
Therefore, each of the four vertices of $H^5_{12}$  of the
type $\widetilde F_4$
belongs to either one  or to at least 16  fundamental simplices.
But 16 fundamental simplices have 16 vertices of the type $\widetilde F_4$.
This is impossible.

\subsubsection{There is no decomposition with symbol
$(H^5_7,H^5_{11}),\ N=6$.}

Suppose that there exists a decomposition with symbol
$(H^5_7,H^5_{11})$ and $N=6$. Then 6 fundamental simplices contain
$3\cdot 6=18$ ideal vertices. These 18 vertices coincide with 5
ideal vertices of $P$.
An ideal vertex of $P$ belongs to at least two fundamental simplices,
since the diagram of $F$ contains no subdiagram $\widetilde D_4$.
No vertex of $P$ belongs to exactly two fundamental simplices,
otherwise by Lemma~\ref{2} $P$ has a non-fundamental dihedral angle.

Suppose that an ideal vertex of $P$ belongs to exactly 3 fundamental
simplices. Then there exists a decomposition of the second type with
$N=3$ and either
$(F,P)=(\widetilde B_4,\widetilde D_4)$ or $(F,P)=(\widetilde
F_4,\widetilde D_4)$.
But $\widetilde D_4$ contains a subdiagram $4A_1$,
and none of the diagrams $\widetilde B_4$ and $\widetilde F_4$
contains this subdiagram.
Moreover, the simplex
$4A_1$ admits no decomposition of the second type with
$N=3$.
Therefore, there is no decomposition with $N=3$ and with the symbols
$(\widetilde B_4,\widetilde D_4)$ and
$(\widetilde F_4,\widetilde D_4)$.
Hence, no ideal vertex of $P$ belongs to exactly 3 fundamental simplices.

Thus, each of five ideal vertices of $P$ belongs to at least four
fundamental simplices. This contradicts to the fact that six fundamental
simplices have 18 ideal vertices.

\subsubsection{There exists a unique decomposition with symbol
$(H^5_4,H^5_{11})$ and $N=20$.}
\label{type2.5}

By Lemma~\ref{similar}, the ideal vertex of $P$
belongs to either exactly one fundamental simplex or
to at least 16 fundamental simplices.
The simplex $P=H^5_{11}$ has 5 ideal vertices.
Since $N=20$, four of the ideal vertices of $P$
belong to exactly one fundamental simplex each,
the rest vertex belongs to 16 fundamental simplices.

\begin{lemma}
There exists at most one decomposition with symbol
$(H^5_4,H^5_{11})$.

\end{lemma}

\begin{proof}
Let $v$ be a vertex of $P$ which belongs to exactly one fundamental
simplex $F_1$. Let $\alpha$ be a facet of $F_1$ opposite to $v$.
The position of $F_1$ with respect to $P$ determines the decomposition.
Hence, it is sufficient to show that the way to put the facet $\alpha$
is unique. Let $\beta_1$,...,$\beta_5$ be the facets of $P$ containing the
vertex $v$.
Clearly, $\alpha$ intersects all but one of
$\beta_1$,...,$\beta_5$
perpendicularly.
The rest facet can be chosen arbitrary in the set
$\{\beta_1$,...,$\beta_5\}$
(the change of this facet leads to a symmetry of the decomposition
and does  not determine  new decompositions).

\end{proof}

\begin{lemma}
There exists a decomposition with symbol
$(H^5_4,H^5_{11}),N=20$.

\end{lemma}

\begin{proof}
The proof is a straightforward calculation in the linear model of
the hyperbolic space.

\end{proof}

{
\begin{itemize}
\item[]

Now we are aimed to describe the decomposition.
The decomposition is very similar to the decomposition
with symbol
$(H^4_3,H^4_9),N=10$.
We can visualize the decomposition as follows.

First, we describe the decomposition in the neighborhood
of the vertex $v$ which is incident to sixteen fundamental simplices
$F_1,...,F_{16}$.
Let $s$ be a horosphere centered in $v$.
Then $s\bigcap P$ is a decomposition with symbol
$(\widetilde D_4,\widetilde D_4), N^4=16$.
The simplex $D_4$ admitting a decomposition can be explained
as
$\{P^4=x_1+x_2+x_3+x_4\le 2,|\quad x_i\ge 0\}$ in Cartesian coordinates.
Then the vertex $(0,0,0,0)$ belongs to twelve fundamental simplices
$f_1,...,f_{12}$.  The decomposition in the neighborhood of $(0,0,0,0)$
coincides with the decomposition with symbol
$(D_3,A_1+A_1+A_1+A_1), N^3=12$.
The rest four fundamental simplices  $f_{13},...,f_{16}$
are cut from $P^4$ by the 3-dimensional faces $x_i=1,
i=1,...,4$.

Let $\bar f_i$ be a facet of $F_i$ ($i=1,...,16$)
opposite to the vertex $v$. Then $\bar f_1$,...,$\bar
f_{12}$ belong to the same hyperplane $\Pi$.
The facets $\bar f_{13},...,\bar f_{16}$ belong to four different
hyperplanes, which intersect  $\Pi$ with dihedral angles
 $\frac{\pi}{3}$.
Reflect the simplices $F_{13},...,F_{16}$ with
respect to $\bar f_{13},...,\bar f_{16}$ respectively
to obtain the rest four fundamental simplices
$F_{17},...,F_{20}$.

\end{itemize}
}

\subsubsection{There exists a unique decomposition with symbol
$(H^8_1,H^8_4)$ and $N=272$.
}
\label{type2.8}

\begin{lemma}
There exists a unique decomposition with symbol $(H^8_1,H^8_4)$ and
$N=272$.

\end{lemma}

\begin{proof}
The existence is proved in~\cite{Ruth2}. Here we reproduce the sckech
of the proof.

Denote the facets of $F=H^8_1$
as shown in Fig.~\ref{H^8_4}.a.
Let $v_i$ be a unit outward normal for the $i$-th  facet of $F$.
Set $v_9=2v_0+v_2+2v_3+3v_4+2v_5+v_6$.
It is easy to check that the vectors $v_1,v_2,...,v_9$
are the outward normals for the simplex $P=H^8_4$ (see
Fig.~\ref{H^8_4}.b).
As it is shown in~\cite{Ruth2}, $v_9$ can be obtained from
$v_2$ by the sequence of the reflections with respect to
the facets orthogonal to
$v_i$, $i=0,...,8$.

\begin{figure}
\label{H^8_4}
\begin{center}
\caption{}
{\small
\setlength{\unitlength}{0.075in}
\begin{picture}(62.15,12.01)
\put(39.61,6.33){{\setbox0=\hbox{3}\lower\ht0\box0}}
\put(36.95,0.00){{\setbox0=\hbox{b)}\lower\ht0\box0}}
\put(0.01,0.00){{\setbox0=\hbox{a)}\lower\ht0\box0}}
\put(54.81,2.00){{\setbox0=\hbox{9}\lower\ht0\box0}}
\put(49.88,1.83){{\setbox0=\hbox{1}\lower\ht0\box0}}
\put(44.68,1.83){{\setbox0=\hbox{2}\lower\ht0\box0}}
\put(39.61,6.33){{\setbox0=\hbox{3}\lower\ht0\box0}}
\put(44.41,11.83){{\setbox0=\hbox{4}\lower\ht0\box0}}
\put(49.35,12.00){{\setbox0=\hbox{5}\lower\ht0\box0}}
\put(54.81,11.83){{\setbox0=\hbox{6}\lower\ht0\box0}}
\put(58.81,8.00){{\setbox0=\hbox{7}\lower\ht0\box0}}
\put(62.15,8.00){{\setbox0=\hbox{8}\lower\ht0\box0}}
\put(61.48,6.33){{\setbox0=\hbox{$\scriptstyle\bullet$}\kern-.4\wd0\lower.5\ht0\box0}}
\put(49.75,9.66){{\setbox0=\hbox{$\scriptstyle\bullet$}\kern-.4\wd0\lower.5\ht0\box0}}
\put(54.81,9.66){{\setbox0=\hbox{$\scriptstyle\bullet$}\kern-.4\wd0\lower.5\ht0\box0}}
\put(58.15,6.33){{\setbox0=\hbox{$\scriptstyle\bullet$}\kern-.4\wd0\lower.5\ht0\box0}}
\put(54.81,3.00){{\setbox0=\hbox{$\scriptstyle\bullet$}\kern-.4\wd0\lower.5\ht0\box0}}
\put(49.88,3.00){{\setbox0=\hbox{$\scriptstyle\bullet$}\kern-.4\wd0\lower.5\ht0\box0}}
\put(44.81,3.00){{\setbox0=\hbox{$\scriptstyle\bullet$}\kern-.4\wd0\lower.5\ht0\box0}}
\put(41.48,6.33){{\setbox0=\hbox{$\scriptstyle\bullet$}\kern-.4\wd0\lower.5\ht0\box0}}
\put(44.81,9.66){{\setbox0=\hbox{$\scriptstyle\bullet$}\kern-.4\wd0\lower.5\ht0\box0}}
\special{em:linewidth 0.014in}
\put(58.15,6.33){\special{em:moveto}}
\put(61.48,6.33){\special{em:lineto}}
\put(41.48,6.33){\special{em:moveto}}
\put(44.81,9.66){\special{em:lineto}}
\put(54.81,9.66){\special{em:lineto}}
\put(58.15,6.33){\special{em:lineto}}
\put(54.81,3.00){\special{em:lineto}}
\put(44.81,3.00){\special{em:lineto}}
\put(41.48,6.33){\special{em:lineto}}
\put(13.21,11.33){{\setbox0=\hbox{0}\lower\ht0\box0}}
\put(25.61,4.66){{\setbox0=\hbox{8}\lower\ht0\box0}}
\put(22.28,4.66){{\setbox0=\hbox{7}\lower\ht0\box0}}
\put(19.08,4.66){{\setbox0=\hbox{6}\lower\ht0\box0}}
\put(15.61,4.66){{\setbox0=\hbox{5}\lower\ht0\box0}}
\put(12.41,4.66){{\setbox0=\hbox{4}\lower\ht0\box0}}
\put(8.95,4.66){{\setbox0=\hbox{3}\lower\ht0\box0}}
\put(5.75,4.66){{\setbox0=\hbox{2}\lower\ht0\box0}}
\put(2.41,4.83){{\setbox0=\hbox{1}\lower\ht0\box0}}
\put(12.41,9.66){{\setbox0=\hbox{$\scriptstyle\bullet$}\kern-.4\wd0\lower.5\ht0\box0}}
\put(12.41,6.33){\special{em:moveto}}
\put(12.41,9.66){\special{em:lineto}}
\put(25.75,6.33){{\setbox0=\hbox{$\scriptstyle\bullet$}\kern-.4\wd0\lower.5\ht0\box0}}
\put(22.41,6.33){{\setbox0=\hbox{$\scriptstyle\bullet$}\kern-.4\wd0\lower.5\ht0\box0}}
\put(19.08,6.33){{\setbox0=\hbox{$\scriptstyle\bullet$}\kern-.4\wd0\lower.5\ht0\box0}}
\put(15.75,6.33){{\setbox0=\hbox{$\scriptstyle\bullet$}\kern-.4\wd0\lower.5\ht0\box0}}
\put(12.41,6.33){{\setbox0=\hbox{$\scriptstyle\bullet$}\kern-.4\wd0\lower.5\ht0\box0}}
\put(9.08,6.33){{\setbox0=\hbox{$\scriptstyle\bullet$}\kern-.4\wd0\lower.5\ht0\box0}}
\put(5.75,6.33){{\setbox0=\hbox{$\scriptstyle\bullet$}\kern-.4\wd0\lower.5\ht0\box0}}
\put(2.41,6.33){{\setbox0=\hbox{$\scriptstyle\bullet$}\kern-.4\wd0\lower.5\ht0\box0}}
\put(2.41,6.33){\special{em:moveto}}
\put(25.75,6.33){\special{em:lineto}}
\end{picture}
}
\end{center}
\end{figure}

To prove the uniquness, consider the Coxeter diagram of
$P=H^8_4$. This diagram has a subdiagram $A_8$ (numbered by $1,2,...,8$
in Fig.~\ref{H^8_4}.b).
The spherical simplex $A_8$ can be decomposed either into one copy of
$A_8$ or into $\frac{2^{14}\cdot 3^5\cdot 5^2\cdot 7}{9!}=2^7\cdot 3\cdot
5>272$ copies of $E_8$. The latter is impossible,
since $N=272$.
Hence, $A_8$ is tiled by one copy of $A_8$
and we have two possibilities:
 the facets of $P$ are numbered either as shown in
in Fig.~\ref{H^8_4}.b or
as shown in
in Fig.~\ref{H^8_4bad}.

\begin{figure}
\label{H^8_4bad}
\begin{center}
\caption{}
{\small
\setlength{\unitlength}{0.075in}
\begin{picture}(23.08,10.01)
\put(15.08,0.00){{\setbox0=\hbox{9}\lower\ht0\box0}}
\put(23.08,5.50){{\setbox0=\hbox{1}\lower\ht0\box0}}
\put(18.81,6.66){{\setbox0=\hbox{2}\lower\ht0\box0}}
\put(14.95,10.00){{\setbox0=\hbox{3}\lower\ht0\box0}}
\put(10.01,10.00){{\setbox0=\hbox{4}\lower\ht0\box0}}
\put(4.28,10.00){{\setbox0=\hbox{5}\lower\ht0\box0}}
\put(0.01,5.16){{\setbox0=\hbox{6}\lower\ht0\box0}}
\put(4.81,0.00){{\setbox0=\hbox{7}\lower\ht0\box0}}
\put(10.15,0.00){{\setbox0=\hbox{8}\lower\ht0\box0}}
\put(22.15,4.33){{\setbox0=\hbox{$\scriptstyle\bullet$}\kern-.4\wd0\lower.5\ht0\box0}}
\put(10.41,7.66){{\setbox0=\hbox{$\scriptstyle\bullet$}\kern-.4\wd0\lower.5\ht0\box0}}
\put(15.48,7.66){{\setbox0=\hbox{$\scriptstyle\bullet$}\kern-.4\wd0\lower.5\ht0\box0}}
\put(18.81,4.33){{\setbox0=\hbox{$\scriptstyle\bullet$}\kern-.4\wd0\lower.5\ht0\box0}}
\put(15.48,1.00){{\setbox0=\hbox{$\scriptstyle\bullet$}\kern-.4\wd0\lower.5\ht0\box0}}
\put(10.55,1.00){{\setbox0=\hbox{$\scriptstyle\bullet$}\kern-.4\wd0\lower.5\ht0\box0}}
\put(5.48,1.00){{\setbox0=\hbox{$\scriptstyle\bullet$}\kern-.4\wd0\lower.5\ht0\box0}}
\put(2.15,4.33){{\setbox0=\hbox{$\scriptstyle\bullet$}\kern-.4\wd0\lower.5\ht0\box0}}
\put(5.48,7.66){{\setbox0=\hbox{$\scriptstyle\bullet$}\kern-.4\wd0\lower.5\ht0\box0}}
\special{em:linewidth 0.014in}
\put(18.81,4.33){\special{em:moveto}}
\put(22.15,4.33){\special{em:lineto}}
\put(2.15,4.33){\special{em:moveto}}
\put(5.48,7.66){\special{em:lineto}}
\put(15.48,7.66){\special{em:lineto}}
\put(18.81,4.33){\special{em:lineto}}
\put(15.48,1.00){\special{em:lineto}}
\put(5.48,1.00){\special{em:lineto}}
\put(2.15,4.33){\special{em:lineto}}
\end{picture}
}
\end{center}
\end{figure}

Suppose that there exists a decomposition corresponding to the letter
numbering. Let $u$ be a unit vector orthogonal to the $9$-th facet.
It is easy to check that $u$ does not belong to the lattice
$\sum\limits_{i=0}^{8}\Z v_i$.
Thus, $u$ can not be obtained form  $v_i$ ($i=0,...,8)$
by a sequence of the reflections with respect to the facets of $F$.
Hence, in this case the decomposition is impossible.

\end{proof}

The proof of the lemma shows that the facets of $P$ are numbered
as it is shown in Fig.~\ref{H^8_4}.
Let $A_i$ be a vertex opposite to the $i$-th facet of $P$.
$P$ has two ideal vertices: $A_3$ and $A_8$.
The fundamental simplex has one ideal vertex.
It is easy to check that $A_3$ belongs to 200
fundamental simplices
and $A_8$ belongs to 72 fundamental simplices.

\subsubsection{There exists a unique decomposition with symbol
$(H^9_1,H^9_3)$ and $N=527$.
}
\label{type2.9}

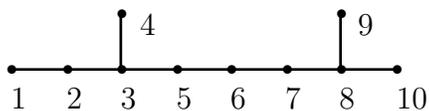
\begin{figure}
\begin{center}
\setlength{\unitlength}{0.0575in}
\begin{picture}(35.08,6.68)
\put(10.01,6.66){{\setbox0=\hbox{$\scriptstyle\bullet$}\kern-.4\wd0\lower.5\ht0\box0}}
\put(11.75,6.66){{\setbox0=\hbox{4}\lower\ht0\box0}}
\special{em:linewidth 0.014in}
\put(10.01,1.66){\special{em:moveto}}
\put(10.01,6.66){\special{em:lineto}}
\put(0.01,0.00){{\setbox0=\hbox{1}\lower\ht0\box0}}
\put(0.01,1.66){{\setbox0=\hbox{$\scriptstyle\bullet$}\kern-.4\wd0\lower.5\ht0\box0}}
\put(5.08,0.00){{\setbox0=\hbox{2}\lower\ht0\box0}}
\put(10.01,0.00){{\setbox0=\hbox{3}\lower\ht0\box0}}
\put(15.08,0.00){{\setbox0=\hbox{5}\lower\ht0\box0}}
\put(20.01,0.00){{\setbox0=\hbox{6}\lower\ht0\box0}}
\put(24.95,0.00){{\setbox0=\hbox{7}\lower\ht0\box0}}
\put(29.88,0.00){{\setbox0=\hbox{8}\lower\ht0\box0}}
\put(35.08,0.00){{\setbox0=\hbox{10}\lower\ht0\box0}}
\put(31.61,6.66){{\setbox0=\hbox{9}\lower\ht0\box0}}
\put(5.08,1.66){{\setbox0=\hbox{$\scriptstyle\bullet$}\kern-.4\wd0\lower.5\ht0\box0}}
\put(10.01,1.66){{\setbox0=\hbox{$\scriptstyle\bullet$}\kern-.4\wd0\lower.5\ht0\box0}}
\put(15.08,1.66){{\setbox0=\hbox{$\scriptstyle\bullet$}\kern-.4\wd0\lower.5\ht0\box0}}
\put(20.01,1.66){{\setbox0=\hbox{$\scriptstyle\bullet$}\kern-.4\wd0\lower.5\ht0\box0}}
\put(25.08,1.66){{\setbox0=\hbox{$\scriptstyle\bullet$}\kern-.4\wd0\lower.5\ht0\box0}}
\put(30.01,1.66){{\setbox0=\hbox{$\scriptstyle\bullet$}\kern-.4\wd0\lower.5\ht0\box0}}
\put(35.08,1.66){{\setbox0=\hbox{$\scriptstyle\bullet$}\kern-.4\wd0\lower.5\ht0\box0}}
\put(30.01,6.66){{\setbox0=\hbox{$\scriptstyle\bullet$}\kern-.4\wd0\lower.5\ht0\box0}}
\put(30.01,6.66){\special{em:moveto}}
\put(30.01,1.66){\special{em:lineto}}
\put(35.08,1.66){\special{em:lineto}}
\put(0.01,1.66){\special{em:moveto}}
\put(30.01,1.66){\special{em:lineto}}
\end{picture}
\end{center}
\caption{The numbering of the facets of the simplex $H^9_3$.}
\label{H^9_3}
\end{figure}

Denote the facets of $P=H^9_3$ as shown in Fig.~\ref{H^9_3}.
Let $v$ be a vertex of $P$ opposite to the fourth facet.
It follows from Subdiagram Property that a fundamental simplex
$F_1$ containing $v$ is unique. There is a unique way to put the simplex
$F_1$ inside $P$. Hence, the decomposition with symbol
$(H^9_1,H^9_3)$ is unique. To check that there exists the decomposition
with this symbol, consider the simplex $F$ in the linear model of $\H^9$.
It is not difficult to check using a computer
that after a number of
reflections we obtain the facet of $P$ opposite to $v$.  It turns out that
the ideal vertex opposite to the first facet belongs to 270 fundamental
simplices, the ideal vertex opposite to the ninth facet belongs to 256
fundamental simplices
and the ideal vertex opposite to the tenth facet belongs to a unique
fundamental simplex.

Thus,  the following Theorem is proved:

\begin{theorem}
A Coxeter decomposition of a hyperbolic simplex of the second type
is uniquely determined by the pair $(F,P)$.
The possibilities for the pairs $(F,P)$ are listed in Table~\ref{pic_hyp2}.

\end{theorem}

\subsection{Decompositions of the third type} \label{3}

Recall that a Coxeter decomposition of the simplex $P$
is called a decomposition of the third type
if $P$ has a non-fundamental dihedral angle
and the decomposition can not be obtained by the inductive algorithm
(see section 2.2).

Let  $\Theta _1(T,P)$ be a decomposition of the simplex $P$
with fundamental simplex $T$
and $\Theta _2(F,T)$ be a decomposition of $T$ with fundamental simplex  $F$.
Clearly, $P$ admits a decomposition $\Theta (F,P)$
with fundamental simplex $F$, where $\Theta (F,P)$ is a superposition
of $\Theta_1(T,P)$ and $\Theta _2(F,T)$.

Denote by $\Theta (F,T)$ the restriction of $\Theta (F,P)$
to the simplex $T$.

\begin{lemma}\label{comp}
Let $\Theta (F,P)$ be a Coxeter decomposition of
$P$ with fundamental simplex $F$.
Suppose that $P$ is decomposed by a mirror of the decomposition
into two simplices $P_1$ and $P_2$.
Suppose that there exists a simplex $T \subseteq P_1$
such that  $\Theta (F,P_1) $ is a superposition of
the  decompositions $\Theta (F,T)$  and
$\Theta _1(T,P_1)$,
where $\Theta _1(T,P_1)$ is a decomposition of $P_1$
with fundamental simplex $T$.
Suppose that the dihedral angles
in the decomposition
$\Theta (F,T)$ are fundamental.
Then the decomposition $\Theta (F,P)$ is a superposition of
$\Theta (F,T)$
and
$\Theta _1(T,P)$,
where $\Theta _1(T,P)$ is a decomposition of $P$ with fundamental simplex
$T$.

\end{lemma}

\begin{proof}
Let
$P=A_0A_1A_2...A_n$, $P_1=A_0A_1...A_{n-1}B$,
$P_2=BA_1A_2...A_{n-1}A_n$ ($B\in A_0A_n$).
Let $\alpha =A_1A_2...A_{n-1}A_n$.
Extend the decomposition $\Theta _1(T,P_1)$ to the decomposition
$\Theta _1(T)$
of the hyperbolic space $\H^n$.
Clearly, any facet of $P$ except $\alpha$
belongs to a mirror of the decomposition $\Theta _1(T)$.
To prove the lemma, it is sufficient to show that $\alpha$
belongs to a mirror of the decomposition
$\Theta _1(T)$.

Let $S$ be a set of mirrors of $\Theta _1(T)$
containing the $(n-2)$-dimensional face
$A_1...A_{n-1}$. Suppose that
$\alpha$ does not belong to $S$.
Then $\alpha$ decomposes a dihedral angle of some fundamental simplex
$T_1$ of the decomposition $\Theta _1(T)$.
This contradicts to the assumption that the dihedral angles
in the decomposition $\Theta (F,T)$ are fundamental.

\end{proof}

\begin{theorem}
Any decomposition of the third type is a superposition
of some non-trivial decompositions.

\end{theorem}

\begin{proof}
Let
$\Theta (F,P)$
be a decomposition of the third type
of a simplex $P$ with fundamental simplex $F$.
Then by the definition of the decomposition of the third type
there exists a mirror decomposing $P$ into two
simplices $P_1$ ¨ $P_2$. Let $\Theta _1$ and  $\Theta _2$ be the
restrictions of $\Theta (F,P)$ to the simplices $P_1$ and $P_2$.  At least
one of $\Theta _1$ and  $\Theta _2$ is a decomposition either of the
second or of the third type (otherwise, $\Theta (F,P)$ could be obtained
by the inductive algorithm).  Suppose that $\Theta _1$ is a decomposition
of the second type.  Then, applying Lemma~\ref{comp}  to the
simplex $T=P_1$, we obtain that  $\Theta (F,P)$  is a superposition of
$\Theta _1$ and $\Theta (P_1,P)$, where $\Theta (P_1,P)$ is a
decomposition of $P$ with fundamental simplex $P_1$.

Suppose that $\Theta _1$  is a decomposition of the third type.
Then $P_1$ is decomposed by a mirror into two smaller simplices.
If the decompositions of these simplices are of the third type,
then these simplices are decomposed into some smaller simplices.
Since
$\Theta (F,P)$ is a decomposition of the third type,
the finite number $k$ of the steps leads to a simplex with the decomposition
of the second type.
Then by
Lemma~\ref{comp}
the decompositions of the simplices
at the step $(k-1)$  is a superposition of some decompositions.
Use Lemma~\ref{comp} again to show that the decomposition
at the step $(k-2)$ is a superposition too.
Applying the Lemma~\ref{comp} for $k$ times,
we obtain the statement of the theorem.

\end{proof}

Thus, any simple decomposition of a hyperbolic simplex
is a decomposition either of the first or of the second type.
Hence, all the simple decompositions are listed in Table~3
and Table~4.

\medskip

\pagebreak

\section{Tables}
\vspace{10pt}
\begin{center}
{\normalsize Table 2:
\footnote{The volumes are rewriten from~\cite{Ruth}.}
Hyperbolic Coxeter simplices in $\H^n$,  $\ge 4$.}\\
\begin{tabular}{|c|c|c|}
\multicolumn{2}{c}{}\\
%
\hline
\begin{tabular}{c}
\phantom{idas}
Coxeter diagram
\phantom{idas}
\\
\end{tabular}
&  Notation &
 Volum  \\
\hline
{\scriptsize
\parbox[c]{3.75cm}%
{
\setlength{\unitlength}{0.083in}
\begin{picture}(17.73,13.45)
\special{em:linewidth 0.0104in}
\put(1.93,0.00){\special{em:moveto}}
\put(2.85,1.83){\special{em:lineto}}
\put(4.78,1.83){\special{em:lineto}}
\put(5.80,0.00){\special{em:lineto}}
\put(6.06,1.76){{\setbox0=\hbox{4}\lower\ht0\box0}}
\put(4.78,1.83){{\setbox0=\hbox{$\scriptstyle\bullet$}\kern-.4\wd0\lower.5\ht0\box0}}
\put(5.80,0.00){{\setbox0=\hbox{$\scriptstyle\bullet$}\kern-.4\wd0\lower.5\ht0\box0}}
\put(3.86,0.00){{\setbox0=\hbox{$\scriptstyle\bullet$}\kern-.4\wd0\lower.5\ht0\box0}}
\put(1.93,0.00){{\setbox0=\hbox{$\scriptstyle\bullet$}\kern-.4\wd0\lower.5\ht0\box0}}
\put(2.85,1.83){{\setbox0=\hbox{$\scriptstyle\bullet$}\kern-.4\wd0\lower.5\ht0\box0}}
\put(1.93,0.00){\special{em:moveto}}
\put(5.80,0.00){\special{em:lineto}}
\put(10.86,4.76){{\setbox0=\hbox{5}\lower\ht0\box0}}
\put(16.73,4.76){{\setbox0=\hbox{5}\lower\ht0\box0}}
\put(17.73,2.93){{\setbox0=\hbox{$\scriptstyle\bullet$}\kern-.4\wd0\lower.5\ht0\box0}}
\put(15.80,2.93){{\setbox0=\hbox{$\scriptstyle\bullet$}\kern-.4\wd0\lower.5\ht0\box0}}
\put(13.86,2.93){{\setbox0=\hbox{$\scriptstyle\bullet$}\kern-.4\wd0\lower.5\ht0\box0}}
\put(11.93,2.93){{\setbox0=\hbox{$\scriptstyle\bullet$}\kern-.4\wd0\lower.5\ht0\box0}}
\put(10.00,2.93){{\setbox0=\hbox{$\scriptstyle\bullet$}\kern-.4\wd0\lower.5\ht0\box0}}
\put(10.00,2.93){\special{em:moveto}}
\put(17.73,2.93){\special{em:lineto}}
\put(6.71,5.20){\special{em:moveto}}
\put(4.86,6.16){\special{em:lineto}}
\put(6.71,7.13){\special{em:lineto}}
\put(1.93,8.10){{\setbox0=\hbox{5}\lower\ht0\box0}}
\put(6.71,7.13){{\setbox0=\hbox{$\scriptstyle\bullet$}\kern-.4\wd0\lower.5\ht0\box0}}
\put(6.71,5.20){{\setbox0=\hbox{$\scriptstyle\bullet$}\kern-.4\wd0\lower.5\ht0\box0}}
\put(4.86,6.16){{\setbox0=\hbox{$\scriptstyle\bullet$}\kern-.4\wd0\lower.5\ht0\box0}}
\put(2.93,6.16){{\setbox0=\hbox{$\scriptstyle\bullet$}\kern-.4\wd0\lower.5\ht0\box0}}
\put(1.00,6.16){{\setbox0=\hbox{$\scriptstyle\bullet$}\kern-.4\wd0\lower.5\ht0\box0}}
\put(1.00,6.16){\special{em:moveto}}
\put(4.86,6.16){\special{em:lineto}}
\put(10.73,11.26){{\setbox0=\hbox{4}\lower\ht0\box0}}
\put(16.60,11.26){{\setbox0=\hbox{5}\lower\ht0\box0}}
\put(17.60,9.40){{\setbox0=\hbox{$\scriptstyle\bullet$}\kern-.4\wd0\lower.5\ht0\box0}}
\put(15.66,9.40){{\setbox0=\hbox{$\scriptstyle\bullet$}\kern-.4\wd0\lower.5\ht0\box0}}
\put(13.73,9.40){{\setbox0=\hbox{$\scriptstyle\bullet$}\kern-.4\wd0\lower.5\ht0\box0}}
\put(11.80,9.40){{\setbox0=\hbox{$\scriptstyle\bullet$}\kern-.4\wd0\lower.5\ht0\box0}}
\put(9.86,9.40){{\setbox0=\hbox{$\scriptstyle\bullet$}\kern-.4\wd0\lower.5\ht0\box0}}
\put(9.86,9.40){\special{em:moveto}}
\put(17.60,9.40){\special{em:lineto}}
\put(6.73,13.43){{\setbox0=\hbox{5}\lower\ht0\box0}}
\put(7.73,11.50){{\setbox0=\hbox{$\scriptstyle\bullet$}\kern-.4\wd0\lower.5\ht0\box0}}
\put(5.80,11.50){{\setbox0=\hbox{$\scriptstyle\bullet$}\kern-.4\wd0\lower.5\ht0\box0}}
\put(3.86,11.50){{\setbox0=\hbox{$\scriptstyle\bullet$}\kern-.4\wd0\lower.5\ht0\box0}}
\put(1.93,11.50){{\setbox0=\hbox{$\scriptstyle\bullet$}\kern-.4\wd0\lower.5\ht0\box0}}
\put(0.00,11.50){{\setbox0=\hbox{$\scriptstyle\bullet$}\kern-.4\wd0\lower.5\ht0\box0}}
\put(0.00,11.50){\special{em:moveto}}
\put(7.73,11.50){\special{em:lineto}}
\end{picture}
\vphantom{$H_a$}}
}
&
\begin{tabular}{c}
\vphantom{$\int\limits^{A^4}$}
$H_1^{(4)}$
\vphantom{$\int\limits^{A^4}$}
\\
$H_2^{(4)}$\\
$H_3^{(4)}$\\
$H_4^{(4)}$\\
\vphantom{$H_{\int\limits_a}$}
$H_5^{(4)}$
\vphantom{$H_{\int\limits_a}$}
\\
\end{tabular}&
\vspace{-10pt}
\begin{tabular}{c}
\vphantom{$\int\limits^{A^4}$}
\vphantom{$H_1^{(4)}$}
$0.00091385226$
\vphantom{$\int\limits^{A^4}$}
\vphantom{$H_1^{(4)}$}
\\
\vphantom{$H_1^{(4)}$}
$0.00776774420$
\vphantom{$H_1^{(4)}$}
\\
\vphantom{$H_1^{(4)}$}
$0.01553548841$
\vphantom{$H_1^{(4)}$}
\\
\vphantom{$H_1^{(4)}$}
$0.02376015874$
\vphantom{$H_1^{(4)}$}
\\
\vphantom{$H_1^{(4)}$}
$0.02513093713$
\vphantom{$H_{\int\limits_a}$}
\end{tabular}
\vspace{8pt}
\\
\hline
\vspace{-3pt}
{\scriptsize
\parbox[c]{3.72cm}{
\setlength{\unitlength}{0.083in}
\begin{picture}(17.65,24.48)
\special{em:linewidth 0.0104in}
\put(3.41,3.38){\special{em:moveto}}
\put(6.01,1.76){\special{em:lineto}}
\put(3.43,0.00){\special{em:lineto}}
\put(3.43,0.00){\special{em:moveto}}
\put(0.60,1.76){\special{em:lineto}}
\put(3.41,3.38){\special{em:lineto}}
\put(3.41,3.38){{\setbox0=\hbox{$\scriptstyle\bullet$}\kern-.4\wd0\lower.5\ht0\box0}}
\put(3.38,0.00){{\setbox0=\hbox{$\scriptstyle\bullet$}\kern-.4\wd0\lower.5\ht0\box0}}
\put(0.81,1.76){{\setbox0=\hbox{$\scriptstyle\bullet$}\kern-.4\wd0\lower.5\ht0\box0}}
\put(3.41,1.76){{\setbox0=\hbox{$\scriptstyle\bullet$}\kern-.4\wd0\lower.5\ht0\box0}}
\put(6.01,1.76){{\setbox0=\hbox{$\scriptstyle\bullet$}\kern-.4\wd0\lower.5\ht0\box0}}
\put(5.90,1.76){\special{em:moveto}}
\put(0.81,1.76){\special{em:lineto}}
\put(10.68,5.13){{\setbox0=\hbox{4}\lower\ht0\box0}}
\put(11.65,3.41){\special{em:moveto}}
\put(12.63,5.35){\special{em:lineto}}
\put(14.65,5.35){\special{em:lineto}}
\put(15.70,3.41){\special{em:lineto}}
\put(16.01,4.96){{\setbox0=\hbox{4}\lower\ht0\box0}}
\put(14.65,5.35){{\setbox0=\hbox{$\scriptstyle\bullet$}\kern-.4\wd0\lower.5\ht0\box0}}
\put(15.70,3.41){{\setbox0=\hbox{$\scriptstyle\bullet$}\kern-.4\wd0\lower.5\ht0\box0}}
\put(13.68,3.41){{\setbox0=\hbox{$\scriptstyle\bullet$}\kern-.4\wd0\lower.5\ht0\box0}}
\put(11.65,3.41){{\setbox0=\hbox{$\scriptstyle\bullet$}\kern-.4\wd0\lower.5\ht0\box0}}
\put(12.63,5.35){{\setbox0=\hbox{$\scriptstyle\bullet$}\kern-.4\wd0\lower.5\ht0\box0}}
\put(11.65,3.41){\special{em:moveto}}
\put(15.70,3.41){\special{em:lineto}}
\put(1.48,8.63){{\setbox0=\hbox{4}\lower\ht0\box0}}
\put(13.81,9.66){\special{em:moveto}}
\put(11.78,8.75){\special{em:lineto}}
\put(15.75,8.65){\special{em:moveto}}
\put(13.81,9.66){\special{em:lineto}}
\put(15.75,10.68){\special{em:lineto}}
\put(12.15,11.96){{\setbox0=\hbox{4}\lower\ht0\box0}}
\put(15.75,10.68){{\setbox0=\hbox{$\scriptstyle\bullet$}\kern-.4\wd0\lower.5\ht0\box0}}
\put(15.75,8.65){{\setbox0=\hbox{$\scriptstyle\bullet$}\kern-.4\wd0\lower.5\ht0\box0}}
\put(13.81,9.66){{\setbox0=\hbox{$\scriptstyle\bullet$}\kern-.4\wd0\lower.5\ht0\box0}}
\put(11.78,8.75){{\setbox0=\hbox{$\scriptstyle\bullet$}\kern-.4\wd0\lower.5\ht0\box0}}
\put(11.70,10.68){{\setbox0=\hbox{$\scriptstyle\bullet$}\kern-.4\wd0\lower.5\ht0\box0}}
\put(11.70,10.68){\special{em:moveto}}
\put(13.81,9.66){\special{em:lineto}}
\put(4.55,15.13){{\setbox0=\hbox{4}\lower\ht0\box0}}
\put(4.40,7.96){\special{em:moveto}}
\put(6.18,6.83){\special{em:lineto}}
\put(4.40,5.93){\special{em:lineto}}
\put(4.40,5.93){\special{em:moveto}}
\put(2.46,6.95){\special{em:lineto}}
\put(4.40,7.96){\special{em:lineto}}
\put(4.40,7.96){{\setbox0=\hbox{$\scriptstyle\bullet$}\kern-.4\wd0\lower.5\ht0\box0}}
\put(4.40,5.93){{\setbox0=\hbox{$\scriptstyle\bullet$}\kern-.4\wd0\lower.5\ht0\box0}}
\put(2.46,6.95){{\setbox0=\hbox{$\scriptstyle\bullet$}\kern-.4\wd0\lower.5\ht0\box0}}
\put(0.43,6.95){{\setbox0=\hbox{$\scriptstyle\bullet$}\kern-.4\wd0\lower.5\ht0\box0}}
\put(6.18,6.83){{\setbox0=\hbox{$\scriptstyle\bullet$}\kern-.4\wd0\lower.5\ht0\box0}}
\put(0.43,6.95){\special{em:moveto}}
\put(2.46,6.95){\special{em:lineto}}
\put(4.45,18.90){\special{em:moveto}}
\put(6.23,17.76){\special{em:lineto}}
\put(4.45,16.86){\special{em:lineto}}
\put(6.00,11.53){\special{em:moveto}}
\put(5.01,12.03){\special{em:lineto}}
\put(4.05,12.55){\special{em:lineto}}
\put(6.00,13.56){\special{em:lineto}}
\put(0.81,14.46){{\setbox0=\hbox{4}\lower\ht0\box0}}
\put(6.00,13.56){{\setbox0=\hbox{$\scriptstyle\bullet$}\kern-.4\wd0\lower.5\ht0\box0}}
\put(6.00,11.53){{\setbox0=\hbox{$\scriptstyle\bullet$}\kern-.4\wd0\lower.5\ht0\box0}}
\put(4.05,12.55){{\setbox0=\hbox{$\scriptstyle\bullet$}\kern-.4\wd0\lower.5\ht0\box0}}
\put(2.03,12.55){{\setbox0=\hbox{$\scriptstyle\bullet$}\kern-.4\wd0\lower.5\ht0\box0}}
\put(0.00,12.55){{\setbox0=\hbox{$\scriptstyle\bullet$}\kern-.4\wd0\lower.5\ht0\box0}}
\put(0.00,12.55){\special{em:moveto}}
\put(4.05,12.55){\special{em:lineto}}
\put(16.26,14.30){\special{em:moveto}}
\put(14.31,15.31){\special{em:lineto}}
\put(16.26,16.33){\special{em:lineto}}
\put(13.08,16.96){{\setbox0=\hbox{4}\lower\ht0\box0}}
\put(16.26,16.33){{\setbox0=\hbox{$\scriptstyle\bullet$}\kern-.4\wd0\lower.5\ht0\box0}}
\put(16.26,14.30){{\setbox0=\hbox{$\scriptstyle\bullet$}\kern-.4\wd0\lower.5\ht0\box0}}
\put(14.31,15.31){{\setbox0=\hbox{$\scriptstyle\bullet$}\kern-.4\wd0\lower.5\ht0\box0}}
\put(12.30,15.31){{\setbox0=\hbox{$\scriptstyle\bullet$}\kern-.4\wd0\lower.5\ht0\box0}}
\put(10.26,15.31){{\setbox0=\hbox{$\scriptstyle\bullet$}\kern-.4\wd0\lower.5\ht0\box0}}
\put(10.26,15.31){\special{em:moveto}}
\put(14.31,15.31){\special{em:lineto}}
\put(4.45,16.86){\special{em:moveto}}
\put(2.51,17.88){\special{em:lineto}}
\put(4.45,18.90){\special{em:lineto}}
\put(4.45,18.90){{\setbox0=\hbox{$\scriptstyle\bullet$}\kern-.4\wd0\lower.5\ht0\box0}}
\put(4.45,16.86){{\setbox0=\hbox{$\scriptstyle\bullet$}\kern-.4\wd0\lower.5\ht0\box0}}
\put(2.51,17.88){{\setbox0=\hbox{$\scriptstyle\bullet$}\kern-.4\wd0\lower.5\ht0\box0}}
\put(0.48,17.88){{\setbox0=\hbox{$\scriptstyle\bullet$}\kern-.4\wd0\lower.5\ht0\box0}}
\put(6.23,17.76){{\setbox0=\hbox{$\scriptstyle\bullet$}\kern-.4\wd0\lower.5\ht0\box0}}
\put(0.48,17.88){\special{em:moveto}}
\put(2.51,17.88){\special{em:lineto}}
\put(6.30,21.50){\special{em:moveto}}
\put(5.31,22.01){\special{em:lineto}}
\put(4.35,22.51){\special{em:lineto}}
\put(6.30,23.53){\special{em:lineto}}
\put(4.28,24.46){{\setbox0=\hbox{4}\lower\ht0\box0}}
\put(6.30,23.53){{\setbox0=\hbox{$\scriptstyle\bullet$}\kern-.4\wd0\lower.5\ht0\box0}}
\put(6.30,21.50){{\setbox0=\hbox{$\scriptstyle\bullet$}\kern-.4\wd0\lower.5\ht0\box0}}
\put(4.35,22.51){{\setbox0=\hbox{$\scriptstyle\bullet$}\kern-.4\wd0\lower.5\ht0\box0}}
\put(2.33,22.51){{\setbox0=\hbox{$\scriptstyle\bullet$}\kern-.4\wd0\lower.5\ht0\box0}}
\put(0.30,22.51){{\setbox0=\hbox{$\scriptstyle\bullet$}\kern-.4\wd0\lower.5\ht0\box0}}
\put(0.30,22.51){\special{em:moveto}}
\put(4.35,22.51){\special{em:lineto}}
\put(12.41,22.13){{\setbox0=\hbox{4}\lower\ht0\box0}}
\put(16.41,22.13){{\setbox0=\hbox{4}\lower\ht0\box0}}
\put(17.65,20.58){{\setbox0=\hbox{$\scriptstyle\bullet$}\kern-.4\wd0\lower.5\ht0\box0}}
\put(15.61,20.58){{\setbox0=\hbox{$\scriptstyle\bullet$}\kern-.4\wd0\lower.5\ht0\box0}}
\put(13.58,20.58){{\setbox0=\hbox{$\scriptstyle\bullet$}\kern-.4\wd0\lower.5\ht0\box0}}
\put(11.55,20.58){{\setbox0=\hbox{$\scriptstyle\bullet$}\kern-.4\wd0\lower.5\ht0\box0}}
\put(9.51,20.58){{\setbox0=\hbox{$\scriptstyle\bullet$}\kern-.4\wd0\lower.5\ht0\box0}}
\put(9.51,20.58){\special{em:moveto}}
\put(17.65,20.58){\special{em:lineto}}
\end{picture}}
\vphantom{$H_{\int\limits_{A_A}}$}
}
&
\begin{tabular}{c}
\vphantom{$\int\limits$}
$H_1^4$
\vphantom{$\int\limits$}
\\
\vphantom{$H_1^{(4)}$}
$H_2^4$
\vphantom{$H_1^{(4)}$}
\\
\vphantom{$H_1^{(4)}$}
$H_3^4$
\vphantom{$H_1^{(4)}$}
\\
\vphantom{$H_1^{(4)}$}
$H_4^4$
\vphantom{$H_1^{(4)}$}
\\
\vphantom{$H_1^{(4)}$}
$H_5^4$
\vphantom{$H_1^{(4)}$}
\\
\vphantom{$H_1^{(4)}$}
$H_6^4$
\vphantom{$H_1^{(4)}$}
\\
\vphantom{$H_1^{(4)}$}
$H_7^4$
\vphantom{$H_1^{(4)}$}
\\
\vphantom{$H_1^{(4)}$}
$H_8^4$
\vphantom{$H_1^{(4)}$}
\\
\vphantom{$H_1^{(4)}$}
$H_9^4$
\vphantom{$H_1^{(4)}$}
\\
\end{tabular}&
\begin{tabular}{c}
\vphantom{$\int\limits^{A}$}
$0.00685389195$
\vphantom{$\int\limits^{A}$}
\\
\vphantom{$H_1^{(4)}$}
$0.01142315324$
\vphantom{$H_1^{(4)}$}
\\
\vphantom{$H_1^{(4)}$}
$0.01370778389$
\vphantom{$H_1^{(4)}$}
\\
\vphantom{$H_1^{(4)}$}
$0.02284630648$
\vphantom{$H_1^{(4)}$}
\\
\vphantom{$H_1^{(4)}$}
$0.03426945973$
\vphantom{$H_1^{(4)}$}
\\
\vphantom{$H_1^{(4)}$}
$0.06853891945$
\vphantom{$H_1^{(4)}$}
\\
\vphantom{$H_1^{(4)}$}
$0.06853891945$
\vphantom{$H_1^{(4)}$}
\\
\vphantom{$H_1^{(4)}$}
$0.09138522594$
\vphantom{$H_1^{(4)}$}
\\
\vphantom{$H_{\int\limits_a}$}
\vphantom{$H_1^{(4)}$}
$0.13707783890$
\vphantom{$H_1^{(4)}$}
\vphantom{$H_{\int\limits_a}$}
\\
\end{tabular}\\
\hline
{\scriptsize
\parbox[c]{4.1cm}{
\setlength{\unitlength}{0.077in}
\begin{picture}(20.63,34.68)
\put(16.16,4.33){{\setbox0=\hbox{4}\lower\ht0\box0}}
\put(16.03,0.00){{\setbox0=\hbox{4}\lower\ht0\box0}}
\put(17.20,2.45){{\setbox0=\hbox{$\scriptstyle\bullet$}\kern-.4\wd0\lower.5\ht0\box0}}
\special{em:linewidth 0.0104in}
\put(15.20,2.45){\special{em:moveto}}
\put(17.20,2.45){\special{em:lineto}}
\put(19.06,1.48){\special{em:lineto}}
\put(17.20,0.51){\special{em:lineto}}
\put(15.20,0.51){\special{em:lineto}}
\put(15.20,0.51){\special{em:moveto}}
\put(13.33,1.48){\special{em:lineto}}
\put(15.20,2.45){\special{em:lineto}}
\put(15.20,2.45){{\setbox0=\hbox{$\scriptstyle\bullet$}\kern-.4\wd0\lower.5\ht0\box0}}
\put(15.20,0.51){{\setbox0=\hbox{$\scriptstyle\bullet$}\kern-.4\wd0\lower.5\ht0\box0}}
\put(13.33,1.48){{\setbox0=\hbox{$\scriptstyle\bullet$}\kern-.4\wd0\lower.5\ht0\box0}}
\put(19.06,1.48){{\setbox0=\hbox{$\scriptstyle\bullet$}\kern-.4\wd0\lower.5\ht0\box0}}
\put(17.20,0.51){{\setbox0=\hbox{$\scriptstyle\bullet$}\kern-.4\wd0\lower.5\ht0\box0}}
\put(6.35,1.46){{\setbox0=\hbox{$\scriptstyle\bullet$}\kern-.4\wd0\lower.5\ht0\box0}}
\put(6.98,4.00){{\setbox0=\hbox{$\scriptstyle\bullet$}\kern-.4\wd0\lower.5\ht0\box0}}
\put(3.45,1.46){{\setbox0=\hbox{$\scriptstyle\bullet$}\kern-.4\wd0\lower.5\ht0\box0}}
\put(2.81,4.00){{\setbox0=\hbox{$\scriptstyle\bullet$}\kern-.4\wd0\lower.5\ht0\box0}}
\put(4.83,4.95){{\setbox0=\hbox{$\scriptstyle\bullet$}\kern-.4\wd0\lower.5\ht0\box0}}
\put(4.83,2.88){\special{em:moveto}}
\put(4.83,4.95){\special{em:lineto}}
\put(3.45,1.46){\special{em:moveto}}
\put(4.83,2.88){\special{em:lineto}}
\put(6.35,1.46){\special{em:lineto}}
\put(2.81,4.00){\special{em:moveto}}
\put(4.83,2.75){\special{em:lineto}}
\put(6.98,4.00){\special{em:lineto}}
\put(4.83,2.88){{\setbox0=\hbox{$\scriptstyle\bullet$}\kern-.4\wd0\lower.5\ht0\box0}}
\put(14.56,8.66){{\setbox0=\hbox{4}\lower\ht0\box0}}
\put(17.66,6.88){\special{em:moveto}}
\put(17.66,5.81){\special{em:lineto}}
\put(17.66,6.88){\special{em:moveto}}
\put(17.66,8.03){\special{em:lineto}}
\put(17.66,8.03){{\setbox0=\hbox{$\scriptstyle\bullet$}\kern-.4\wd0\lower.5\ht0\box0}}
\put(17.66,5.81){{\setbox0=\hbox{$\scriptstyle\bullet$}\kern-.4\wd0\lower.5\ht0\box0}}
\put(19.60,6.88){{\setbox0=\hbox{$\scriptstyle\bullet$}\kern-.4\wd0\lower.5\ht0\box0}}
\put(17.66,6.88){{\setbox0=\hbox{$\scriptstyle\bullet$}\kern-.4\wd0\lower.5\ht0\box0}}
\put(15.73,6.88){{\setbox0=\hbox{$\scriptstyle\bullet$}\kern-.4\wd0\lower.5\ht0\box0}}
\put(13.80,6.88){{\setbox0=\hbox{$\scriptstyle\bullet$}\kern-.4\wd0\lower.5\ht0\box0}}
\put(13.80,6.88){\special{em:moveto}}
\put(19.60,6.88){\special{em:lineto}}
\put(4.83,11.66){{\setbox0=\hbox{4}\lower\ht0\box0}}
\put(5.86,9.93){{\setbox0=\hbox{$\scriptstyle\bullet$}\kern-.4\wd0\lower.5\ht0\box0}}
\put(3.86,9.93){\special{em:moveto}}
\put(5.86,9.93){\special{em:lineto}}
\put(7.73,8.96){\special{em:lineto}}
\put(5.86,8.00){\special{em:lineto}}
\put(3.86,8.00){\special{em:lineto}}
\put(3.86,8.00){\special{em:moveto}}
\put(2.00,8.96){\special{em:lineto}}
\put(3.86,9.93){\special{em:lineto}}
\put(3.86,9.93){{\setbox0=\hbox{$\scriptstyle\bullet$}\kern-.4\wd0\lower.5\ht0\box0}}
\put(3.86,8.00){{\setbox0=\hbox{$\scriptstyle\bullet$}\kern-.4\wd0\lower.5\ht0\box0}}
\put(2.00,8.96){{\setbox0=\hbox{$\scriptstyle\bullet$}\kern-.4\wd0\lower.5\ht0\box0}}
\put(7.73,8.96){{\setbox0=\hbox{$\scriptstyle\bullet$}\kern-.4\wd0\lower.5\ht0\box0}}
\put(5.86,8.00){{\setbox0=\hbox{$\scriptstyle\bullet$}\kern-.4\wd0\lower.5\ht0\box0}}
\put(13.76,14.33){{\setbox0=\hbox{4}\lower\ht0\box0}}
\put(16.70,12.50){\special{em:moveto}}
\put(16.70,13.66){\special{em:lineto}}
\put(16.70,13.66){{\setbox0=\hbox{$\scriptstyle\bullet$}\kern-.4\wd0\lower.5\ht0\box0}}
\put(19.76,14.33){{\setbox0=\hbox{4}\lower\ht0\box0}}
\put(20.63,12.50){{\setbox0=\hbox{$\scriptstyle\bullet$}\kern-.4\wd0\lower.5\ht0\box0}}
\put(18.63,12.50){{\setbox0=\hbox{$\scriptstyle\bullet$}\kern-.4\wd0\lower.5\ht0\box0}}
\put(16.70,12.50){{\setbox0=\hbox{$\scriptstyle\bullet$}\kern-.4\wd0\lower.5\ht0\box0}}
\put(14.76,12.50){{\setbox0=\hbox{$\scriptstyle\bullet$}\kern-.4\wd0\lower.5\ht0\box0}}
\put(12.83,12.50){{\setbox0=\hbox{$\scriptstyle\bullet$}\kern-.4\wd0\lower.5\ht0\box0}}
\put(12.83,12.50){\special{em:moveto}}
\put(20.63,12.50){\special{em:lineto}}
\put(3.20,16.13){{\setbox0=\hbox{$\scriptstyle\bullet$}\kern-.4\wd0\lower.5\ht0\box0}}
\put(9.00,15.16){\special{em:moveto}}
\put(7.13,16.13){\special{em:lineto}}
\put(9.00,17.10){\special{em:lineto}}
\put(4.30,18.16){{\setbox0=\hbox{4}\lower\ht0\box0}}
\put(9.00,17.10){{\setbox0=\hbox{$\scriptstyle\bullet$}\kern-.4\wd0\lower.5\ht0\box0}}
\put(9.00,15.16){{\setbox0=\hbox{$\scriptstyle\bullet$}\kern-.4\wd0\lower.5\ht0\box0}}
\put(7.13,16.13){{\setbox0=\hbox{$\scriptstyle\bullet$}\kern-.4\wd0\lower.5\ht0\box0}}
\put(5.20,16.13){{\setbox0=\hbox{$\scriptstyle\bullet$}\kern-.4\wd0\lower.5\ht0\box0}}
\put(1.26,16.13){{\setbox0=\hbox{$\scriptstyle\bullet$}\kern-.4\wd0\lower.5\ht0\box0}}
\put(1.26,16.13){\special{em:moveto}}
\put(7.13,16.13){\special{em:lineto}}
\put(19.13,19.61){{\setbox0=\hbox{$\scriptstyle\bullet$}\kern-.4\wd0\lower.5\ht0\box0}}
\put(17.11,19.61){\special{em:moveto}}
\put(19.13,19.61){\special{em:lineto}}
\put(19.13,17.68){\special{em:lineto}}
\put(17.11,17.68){\special{em:lineto}}
\put(17.11,17.68){\special{em:moveto}}
\put(15.26,18.65){\special{em:lineto}}
\put(17.11,19.61){\special{em:lineto}}
\put(17.11,19.61){{\setbox0=\hbox{$\scriptstyle\bullet$}\kern-.4\wd0\lower.5\ht0\box0}}
\put(17.11,17.68){{\setbox0=\hbox{$\scriptstyle\bullet$}\kern-.4\wd0\lower.5\ht0\box0}}
\put(15.26,18.65){{\setbox0=\hbox{$\scriptstyle\bullet$}\kern-.4\wd0\lower.5\ht0\box0}}
\put(13.33,18.65){{\setbox0=\hbox{$\scriptstyle\bullet$}\kern-.4\wd0\lower.5\ht0\box0}}
\put(19.13,17.68){{\setbox0=\hbox{$\scriptstyle\bullet$}\kern-.4\wd0\lower.5\ht0\box0}}
\put(13.33,18.65){\special{em:moveto}}
\put(15.26,18.65){\special{em:lineto}}
\put(9.73,22.13){{\setbox0=\hbox{$\scriptstyle\bullet$}\kern-.4\wd0\lower.5\ht0\box0}}
\put(2.70,24.00){{\setbox0=\hbox{4}\lower\ht0\box0}}
\put(7.80,22.13){{\setbox0=\hbox{$\scriptstyle\bullet$}\kern-.4\wd0\lower.5\ht0\box0}}
\put(5.80,22.13){{\setbox0=\hbox{$\scriptstyle\bullet$}\kern-.4\wd0\lower.5\ht0\box0}}
\put(3.86,22.13){{\setbox0=\hbox{$\scriptstyle\bullet$}\kern-.4\wd0\lower.5\ht0\box0}}
\put(1.93,22.13){{\setbox0=\hbox{$\scriptstyle\bullet$}\kern-.4\wd0\lower.5\ht0\box0}}
\put(0.00,22.13){{\setbox0=\hbox{$\scriptstyle\bullet$}\kern-.4\wd0\lower.5\ht0\box0}}
\put(0.00,22.13){\special{em:moveto}}
\put(9.73,22.13){\special{em:lineto}}
\put(17.40,25.33){\special{em:moveto}}
\put(17.40,24.26){\special{em:lineto}}
\put(17.40,25.33){\special{em:moveto}}
\put(17.40,26.50){\special{em:lineto}}
\put(17.40,26.50){{\setbox0=\hbox{$\scriptstyle\bullet$}\kern-.4\wd0\lower.5\ht0\box0}}
\put(8.56,24.16){{\setbox0=\hbox{4}\lower\ht0\box0}}
\put(17.40,24.26){{\setbox0=\hbox{$\scriptstyle\bullet$}\kern-.4\wd0\lower.5\ht0\box0}}
\put(19.33,25.33){{\setbox0=\hbox{$\scriptstyle\bullet$}\kern-.4\wd0\lower.5\ht0\box0}}
\put(17.40,25.33){{\setbox0=\hbox{$\scriptstyle\bullet$}\kern-.4\wd0\lower.5\ht0\box0}}
\put(15.46,25.33){{\setbox0=\hbox{$\scriptstyle\bullet$}\kern-.4\wd0\lower.5\ht0\box0}}
\put(13.53,25.33){{\setbox0=\hbox{$\scriptstyle\bullet$}\kern-.4\wd0\lower.5\ht0\box0}}
\put(13.53,25.33){\special{em:moveto}}
\put(19.33,25.33){\special{em:lineto}}
\put(16.66,30.41){\special{em:moveto}}
\put(16.66,31.58){\special{em:lineto}}
\put(9.76,27.15){{\setbox0=\hbox{$\scriptstyle\bullet$}\kern-.4\wd0\lower.5\ht0\box0}}
\put(4.56,29.33){{\setbox0=\hbox{4}\lower\ht0\box0}}
\put(7.83,27.15){{\setbox0=\hbox{$\scriptstyle\bullet$}\kern-.4\wd0\lower.5\ht0\box0}}
\put(5.81,27.15){{\setbox0=\hbox{$\scriptstyle\bullet$}\kern-.4\wd0\lower.5\ht0\box0}}
\put(3.88,27.15){{\setbox0=\hbox{$\scriptstyle\bullet$}\kern-.4\wd0\lower.5\ht0\box0}}
\put(1.95,27.15){{\setbox0=\hbox{$\scriptstyle\bullet$}\kern-.4\wd0\lower.5\ht0\box0}}
\put(0.01,27.15){{\setbox0=\hbox{$\scriptstyle\bullet$}\kern-.4\wd0\lower.5\ht0\box0}}
\put(0.01,27.15){\special{em:moveto}}
\put(9.76,27.15){\special{em:lineto}}
\put(16.66,31.58){{\setbox0=\hbox{$\scriptstyle\bullet$}\kern-.4\wd0\lower.5\ht0\box0}}
\put(19.23,32.33){{\setbox0=\hbox{4}\lower\ht0\box0}}
\put(20.60,30.41){{\setbox0=\hbox{$\scriptstyle\bullet$}\kern-.4\wd0\lower.5\ht0\box0}}
\put(18.60,30.41){{\setbox0=\hbox{$\scriptstyle\bullet$}\kern-.4\wd0\lower.5\ht0\box0}}
\put(16.66,30.41){{\setbox0=\hbox{$\scriptstyle\bullet$}\kern-.4\wd0\lower.5\ht0\box0}}
\put(14.73,30.41){{\setbox0=\hbox{$\scriptstyle\bullet$}\kern-.4\wd0\lower.5\ht0\box0}}
\put(12.80,30.41){{\setbox0=\hbox{$\scriptstyle\bullet$}\kern-.4\wd0\lower.5\ht0\box0}}
\put(12.80,30.41){\special{em:moveto}}
\put(20.60,30.41){\special{em:lineto}}
\put(9.73,32.68){{\setbox0=\hbox{$\scriptstyle\bullet$}\kern-.4\wd0\lower.5\ht0\box0}}
\put(6.70,34.66){{\setbox0=\hbox{4}\lower\ht0\box0}}
\put(7.80,32.68){{\setbox0=\hbox{$\scriptstyle\bullet$}\kern-.4\wd0\lower.5\ht0\box0}}
\put(5.80,32.68){{\setbox0=\hbox{$\scriptstyle\bullet$}\kern-.4\wd0\lower.5\ht0\box0}}
\put(3.86,32.68){{\setbox0=\hbox{$\scriptstyle\bullet$}\kern-.4\wd0\lower.5\ht0\box0}}
\put(1.93,32.68){{\setbox0=\hbox{$\scriptstyle\bullet$}\kern-.4\wd0\lower.5\ht0\box0}}
\put(0.00,32.68){{\setbox0=\hbox{$\scriptstyle\bullet$}\kern-.4\wd0\lower.5\ht0\box0}}
\put(0.00,32.68){\special{em:moveto}}
\put(9.73,32.68){\special{em:lineto}}
\end{picture}}
\vphantom{$H_{\int\limits_{A_A}}$}
}
&
\begin{tabular}{c}
\vphantom{$\int\limits^{A}$}
$H_1^5$
\vphantom{$\int\limits^{A}$}
\\
\vphantom{$H_1^{(4)}$}
$H_2^5$
\vphantom{$H_1^{(4)}$}
\\
\vphantom{$H_1^{(4)}$}
$H_3^5$
\vphantom{$H_1^{(4)}$}
\\
\vphantom{$H_1^{(4)}$}
$H_4^5$
\vphantom{$H_1^{(4)}$}
\\
\vphantom{$H_1^{(4)}$}
$H_5^5$
\vphantom{$H_1^{(4)}$}
\\
\vphantom{$H_1^{(4)}$}
$H_6^5$
\vphantom{$H_1^{(4)}$}
\\
\vphantom{$H_1^{(4)}$}
$H_7^5$
\vphantom{$H_1^{(4)}$}
\\
\vphantom{$H_1^{(4)}$}
$H_8^5$
\vphantom{$H_1^{(4)}$}
\\
\vphantom{$H_1^{(4)}$}
$H_9^5$
\vphantom{$H_1^{(4)}$}
\\
\vphantom{$H_1^{(4)}$}
$H_{10}^5$
\vphantom{$H_1^{(4)}$}
\\
\vphantom{$H_1^{(4)}$}
$H_{11}^5$
\vphantom{$H_1^{(4)}$}
\\
\vphantom{$H_{\int\limits_{A_A}}^{{A^A}^i}$}
$H_{12}^5$
\vphantom{$H_{\int\limits_{A_A}}^{{A^A}^i}$}
\\
\end{tabular}&
\begin{tabular}{c}
\vphantom{$\int\limits^{A}$}
$0.0001826041$
\vphantom{$\int\limits^{A}$}
\\
\vphantom{$H_1^{(4)}$}
$0.0005478123$
\vphantom{$H_1^{(4)}$}
\\
\vphantom{$H_1^{(4)}$}
$0.0009130206$
\vphantom{$H_1^{(4)}$}
\\
\vphantom{$H_1^{(4)}$}
$0.0010956247$
\vphantom{$H_1^{(4)}$}
\\
\vphantom{$H_1^{(4)}$}
$0.0018260413$
\vphantom{$H_1^{(4)}$}
\\
\vphantom{$H_1^{(4)}$}
$0.0020740519$
\vphantom{$H_1^{(4)}$}
\\
\vphantom{$H_1^{(4)}$}
$0.0036520826$
\vphantom{$H_1^{(4)}$}
\\
\vphantom{$H_1^{(4)}$}
$0.0054781239$
\vphantom{$H_1^{(4)}$}
\\
\vphantom{$H_1^{(4)}$}
$0.0075726186$
\vphantom{$H_1^{(4)}$}
\\
\vphantom{$H_1^{(4)}$}
$0.0109562478$
\vphantom{$H_1^{(4)}$}
\\
\vphantom{$H_1^{(4)}$}
$0.0219124956$
\vphantom{$H_1^{(4)}$}
\\
\vphantom{$H_{\int\limits_{A_A}}^{{A^A}^i}$}
$0.0292166608$
\vphantom{$H_{\int\limits_{A_A}}^{{A^A}^i}$}
\\
\end{tabular}
\\
\hline
\end{tabular}
\begin{tabular}{|c|c|c|}
\multicolumn{3}{c}{{\normalsize Table 2, Continued.}}\\
\multicolumn{2}{c}{}\\
%
\hline
\begin{tabular}{c}
\phantom{asdas}
 Coxeter diagram
\phantom{asdas}
\\
\end{tabular}
&  Notation &
 Volum  \\
\hline
{\scriptsize
\parbox[c]{3cm}{
\setlength{\unitlength}{0.06in}
\begin{picture}(19.78,7.36)
\put(0.61,1.08){{\setbox0=\hbox{$\scriptstyle\bullet$}\kern-.4\wd0\lower.5\ht0\box0}}
\special{em:linewidth 0.0104in}
\put(2.70,1.08){\special{em:moveto}}
\put(0.61,1.08){\special{em:lineto}}
\put(17.53,4.48){{\setbox0=\hbox{$\scriptstyle\bullet$}\kern-.4\wd0\lower.5\ht0\box0}}
\put(17.53,3.18){\special{em:moveto}}
\put(17.53,4.48){\special{em:lineto}}
\put(8.76,5.51){{\setbox0=\hbox{$\scriptstyle\bullet$}\kern-.4\wd0\lower.5\ht0\box0}}
\put(4.33,5.51){\special{em:moveto}}
\put(4.33,6.81){\special{em:lineto}}
\put(4.33,6.81){{\setbox0=\hbox{$\scriptstyle\bullet$}\kern-.4\wd0\lower.5\ht0\box0}}
\put(9.65,7.35){{\setbox0=\hbox{4}\lower\ht0\box0}}
\put(10.83,5.51){{\setbox0=\hbox{$\scriptstyle\bullet$}\kern-.4\wd0\lower.5\ht0\box0}}
\put(6.50,5.51){{\setbox0=\hbox{$\scriptstyle\bullet$}\kern-.4\wd0\lower.5\ht0\box0}}
\put(4.33,5.51){{\setbox0=\hbox{$\scriptstyle\bullet$}\kern-.4\wd0\lower.5\ht0\box0}}
\put(2.16,5.51){{\setbox0=\hbox{$\scriptstyle\bullet$}\kern-.4\wd0\lower.5\ht0\box0}}
\put(0.00,5.51){{\setbox0=\hbox{$\scriptstyle\bullet$}\kern-.4\wd0\lower.5\ht0\box0}}
\put(0.00,5.51){\special{em:moveto}}
\put(10.83,5.51){\special{em:lineto}}
\put(6.93,2.16){{\setbox0=\hbox{$\scriptstyle\bullet$}\kern-.4\wd0\lower.5\ht0\box0}}
\put(4.68,2.16){\special{em:moveto}}
\put(6.93,2.16){\special{em:lineto}}
\put(9.01,1.08){\special{em:lineto}}
\put(6.93,0.00){\special{em:lineto}}
\put(4.68,0.00){\special{em:lineto}}
\put(4.68,0.00){\special{em:moveto}}
\put(2.60,1.08){\special{em:lineto}}
\put(4.68,2.16){\special{em:lineto}}
\put(4.68,2.16){{\setbox0=\hbox{$\scriptstyle\bullet$}\kern-.4\wd0\lower.5\ht0\box0}}
\put(4.68,0.00){{\setbox0=\hbox{$\scriptstyle\bullet$}\kern-.4\wd0\lower.5\ht0\box0}}
\put(2.60,1.08){{\setbox0=\hbox{$\scriptstyle\bullet$}\kern-.4\wd0\lower.5\ht0\box0}}
\put(9.01,1.08){{\setbox0=\hbox{$\scriptstyle\bullet$}\kern-.4\wd0\lower.5\ht0\box0}}
\put(6.93,0.00){{\setbox0=\hbox{$\scriptstyle\bullet$}\kern-.4\wd0\lower.5\ht0\box0}}
\put(15.36,3.18){\special{em:moveto}}
\put(15.36,4.48){\special{em:lineto}}
\put(15.36,4.48){{\setbox0=\hbox{$\scriptstyle\bullet$}\kern-.4\wd0\lower.5\ht0\box0}}
\put(19.78,3.18){{\setbox0=\hbox{$\scriptstyle\bullet$}\kern-.4\wd0\lower.5\ht0\box0}}
\put(17.53,3.18){{\setbox0=\hbox{$\scriptstyle\bullet$}\kern-.4\wd0\lower.5\ht0\box0}}
\put(15.36,3.18){{\setbox0=\hbox{$\scriptstyle\bullet$}\kern-.4\wd0\lower.5\ht0\box0}}
\put(13.20,3.18){{\setbox0=\hbox{$\scriptstyle\bullet$}\kern-.4\wd0\lower.5\ht0\box0}}
\put(11.03,3.18){{\setbox0=\hbox{$\scriptstyle\bullet$}\kern-.4\wd0\lower.5\ht0\box0}}
\put(11.03,3.18){\special{em:moveto}}
\put(19.78,3.18){\special{em:lineto}}
\end{picture}}    }
&
\begin{tabular}{c}
\\
$H_1^6$\\
$H_2^6$\\
$H_3^6$\\
\\
\end{tabular}&

\begin{tabular}{c}
\\
$0.3987432701\times10^{-4}$\\
$0.7974865401\times10^{-4}$\\
$2.9620928633\times10^{-4}$\\
\\
\end{tabular}\\
\hline

{\scriptsize
\parbox[c]{4.22cm}{
\setlength{\unitlength}{0.065in}
\begin{picture}(25.66,9.71)
\put(21.65,0.00){{\setbox0=\hbox{$\scriptstyle\bullet$}\kern-.4\wd0\lower.5\ht0\box0}}
\put(2.36,3.33){{\setbox0=\hbox{$\scriptstyle\bullet$}\kern-.4\wd0\lower.5\ht0\box0}}
\put(23.30,5.36){{\setbox0=\hbox{$\scriptstyle\bullet$}\kern-.4\wd0\lower.5\ht0\box0}}
\put(12.00,1.18){{\setbox0=\hbox{$\scriptstyle\bullet$}\kern-.4\wd0\lower.5\ht0\box0}}
\special{em:linewidth 0.0104in}
\put(14.26,1.18){\special{em:moveto}}
\put(12.00,1.18){\special{em:lineto}}
\put(9.46,4.75){{\setbox0=\hbox{$\scriptstyle\bullet$}\kern-.4\wd0\lower.5\ht0\box0}}
\put(9.46,3.33){\special{em:moveto}}
\put(9.46,4.75){\special{em:lineto}}
\put(21.03,5.36){{\setbox0=\hbox{$\scriptstyle\bullet$}\kern-.4\wd0\lower.5\ht0\box0}}
\put(16.20,5.36){\special{em:moveto}}
\put(16.20,6.78){\special{em:lineto}}
\put(16.20,6.78){{\setbox0=\hbox{$\scriptstyle\bullet$}\kern-.4\wd0\lower.5\ht0\box0}}
\put(24.40,7.41){{\setbox0=\hbox{4}\lower\ht0\box0}}
\put(25.66,5.36){{\setbox0=\hbox{$\scriptstyle\bullet$}\kern-.4\wd0\lower.5\ht0\box0}}
\put(18.56,5.36){{\setbox0=\hbox{$\scriptstyle\bullet$}\kern-.4\wd0\lower.5\ht0\box0}}
\put(16.20,5.36){{\setbox0=\hbox{$\scriptstyle\bullet$}\kern-.4\wd0\lower.5\ht0\box0}}
\put(13.83,5.36){{\setbox0=\hbox{$\scriptstyle\bullet$}\kern-.4\wd0\lower.5\ht0\box0}}
\put(11.46,5.36){{\setbox0=\hbox{$\scriptstyle\bullet$}\kern-.4\wd0\lower.5\ht0\box0}}
\put(11.46,5.36){\special{em:moveto}}
\put(25.66,5.36){\special{em:lineto}}
\put(18.91,2.36){{\setbox0=\hbox{$\scriptstyle\bullet$}\kern-.4\wd0\lower.5\ht0\box0}}
\put(16.45,2.36){\special{em:moveto}}
\put(18.91,2.36){\special{em:lineto}}
\put(21.65,2.36){\special{em:lineto}}
\put(21.65,0.00){\special{em:lineto}}
\put(18.91,0.00){\special{em:lineto}}
\put(16.45,0.00){\special{em:lineto}}
\put(16.45,0.00){\special{em:moveto}}
\put(14.18,1.18){\special{em:lineto}}
\put(16.45,2.36){\special{em:lineto}}
\put(16.45,2.36){{\setbox0=\hbox{$\scriptstyle\bullet$}\kern-.4\wd0\lower.5\ht0\box0}}
\put(16.45,0.00){{\setbox0=\hbox{$\scriptstyle\bullet$}\kern-.4\wd0\lower.5\ht0\box0}}
\put(14.18,1.18){{\setbox0=\hbox{$\scriptstyle\bullet$}\kern-.4\wd0\lower.5\ht0\box0}}
\put(21.65,2.36){{\setbox0=\hbox{$\scriptstyle\bullet$}\kern-.4\wd0\lower.5\ht0\box0}}
\put(18.91,0.00){{\setbox0=\hbox{$\scriptstyle\bullet$}\kern-.4\wd0\lower.5\ht0\box0}}
\put(4.73,3.33){\special{em:moveto}}
\put(4.73,4.86){\special{em:lineto}}
\put(4.73,4.86){{\setbox0=\hbox{$\scriptstyle\bullet$}\kern-.4\wd0\lower.5\ht0\box0}}
\put(11.93,3.33){{\setbox0=\hbox{$\scriptstyle\bullet$}\kern-.4\wd0\lower.5\ht0\box0}}
\put(9.46,3.33){{\setbox0=\hbox{$\scriptstyle\bullet$}\kern-.4\wd0\lower.5\ht0\box0}}
\put(7.10,3.33){{\setbox0=\hbox{$\scriptstyle\bullet$}\kern-.4\wd0\lower.5\ht0\box0}}
\put(4.73,3.33){{\setbox0=\hbox{$\scriptstyle\bullet$}\kern-.4\wd0\lower.5\ht0\box0}}
\put(0.00,3.33){{\setbox0=\hbox{$\scriptstyle\bullet$}\kern-.4\wd0\lower.5\ht0\box0}}
\put(0.00,3.33){\special{em:moveto}}
\put(11.93,3.33){\special{em:lineto}}
\put(11.83,7.36){{\setbox0=\hbox{$\scriptstyle\bullet$}\kern-.4\wd0\lower.5\ht0\box0}}
\put(11.83,9.71){{\setbox0=\hbox{$\scriptstyle\bullet$}\kern-.4\wd0\lower.5\ht0\box0}}
\put(2.56,8.55){{\setbox0=\hbox{$\scriptstyle\bullet$}\kern-.4\wd0\lower.5\ht0\box0}}
\put(11.83,7.36){\special{em:moveto}}
\put(9.66,7.36){\special{em:lineto}}
\put(7.38,8.55){\special{em:lineto}}
\put(9.66,9.71){\special{em:lineto}}
\put(11.83,9.71){\special{em:lineto}}
\put(9.66,9.71){{\setbox0=\hbox{$\scriptstyle\bullet$}\kern-.4\wd0\lower.5\ht0\box0}}
\put(9.66,7.36){{\setbox0=\hbox{$\scriptstyle\bullet$}\kern-.4\wd0\lower.5\ht0\box0}}
\put(7.38,8.55){{\setbox0=\hbox{$\scriptstyle\bullet$}\kern-.4\wd0\lower.5\ht0\box0}}
\put(5.01,8.55){{\setbox0=\hbox{$\scriptstyle\bullet$}\kern-.4\wd0\lower.5\ht0\box0}}
\put(0.20,8.55){{\setbox0=\hbox{$\scriptstyle\bullet$}\kern-.4\wd0\lower.5\ht0\box0}}
\put(0.20,8.55){\special{em:moveto}}
\put(7.38,8.55){\special{em:lineto}}
\end{picture}}    }
&
\begin{tabular}{c}
\\
$H_1^7$\\
$H_2^7$\\
$H_3^7$\\
$H_4^7$\\
\\
\end{tabular}&

\begin{tabular}{c}
$0.1892871372\times10^{-5}$\\
$0.2725266071\times10^{-5}$\\
$0.5450532141\times10^{-5}$\\
$4.1106779054\times10^{-5}$\\
\end{tabular}\\

\hline
{\scriptsize
\parbox[c]{4.77cm}{
\setlength{\unitlength}{0.07in}
\begin{picture}(26.83,8.10)
\put(23.95,1.00){{\setbox0=\hbox{$\scriptstyle\bullet$}\kern-.4\wd0\lower.5\ht0\box0}}
\put(12.25,3.10){{\setbox0=\hbox{$\scriptstyle\bullet$}\kern-.4\wd0\lower.5\ht0\box0}}
\put(26.83,5.10){{\setbox0=\hbox{$\scriptstyle\bullet$}\kern-.4\wd0\lower.5\ht0\box0}}
\put(14.16,7.10){{\setbox0=\hbox{$\scriptstyle\bullet$}\kern-.4\wd0\lower.5\ht0\box0}}
\put(6.08,8.10){{\setbox0=\hbox{$\scriptstyle\bullet$}\kern-.4\wd0\lower.5\ht0\box0}}
\special{em:linewidth 0.0104in}
\put(6.08,7.10){\special{em:moveto}}
\put(6.08,8.10){\special{em:lineto}}
\put(21.95,0.00){{\setbox0=\hbox{$\scriptstyle\bullet$}\kern-.4\wd0\lower.5\ht0\box0}}
\put(2.16,3.10){{\setbox0=\hbox{$\scriptstyle\bullet$}\kern-.4\wd0\lower.5\ht0\box0}}
\put(22.83,5.10){{\setbox0=\hbox{$\scriptstyle\bullet$}\kern-.4\wd0\lower.5\ht0\box0}}
\put(13.80,1.00){{\setbox0=\hbox{$\scriptstyle\bullet$}\kern-.4\wd0\lower.5\ht0\box0}}
\put(15.71,1.00){\special{em:moveto}}
\put(13.80,1.00){\special{em:lineto}}
\put(10.25,4.30){{\setbox0=\hbox{$\scriptstyle\bullet$}\kern-.4\wd0\lower.5\ht0\box0}}
\put(10.25,3.10){\special{em:moveto}}
\put(10.25,4.30){\special{em:lineto}}
\put(20.91,5.10){{\setbox0=\hbox{$\scriptstyle\bullet$}\kern-.4\wd0\lower.5\ht0\box0}}
\put(16.83,5.10){\special{em:moveto}}
\put(16.83,6.30){\special{em:lineto}}
\put(16.83,6.30){{\setbox0=\hbox{$\scriptstyle\bullet$}\kern-.4\wd0\lower.5\ht0\box0}}
\put(25.71,6.90){{\setbox0=\hbox{4}\lower\ht0\box0}}
\put(24.83,5.10){{\setbox0=\hbox{$\scriptstyle\bullet$}\kern-.4\wd0\lower.5\ht0\box0}}
\put(18.83,5.10){{\setbox0=\hbox{$\scriptstyle\bullet$}\kern-.4\wd0\lower.5\ht0\box0}}
\put(16.83,5.10){{\setbox0=\hbox{$\scriptstyle\bullet$}\kern-.4\wd0\lower.5\ht0\box0}}
\put(14.83,5.10){{\setbox0=\hbox{$\scriptstyle\bullet$}\kern-.4\wd0\lower.5\ht0\box0}}
\put(12.83,5.10){{\setbox0=\hbox{$\scriptstyle\bullet$}\kern-.4\wd0\lower.5\ht0\box0}}
\put(12.83,5.10){\special{em:moveto}}
\put(26.83,5.10){\special{em:lineto}}
\put(19.63,2.00){{\setbox0=\hbox{$\scriptstyle\bullet$}\kern-.4\wd0\lower.5\ht0\box0}}
\put(17.55,2.00){\special{em:moveto}}
\put(19.63,2.00){\special{em:lineto}}
\put(21.95,2.00){\special{em:lineto}}
\put(23.95,1.00){\special{em:lineto}}
\put(21.95,0.00){\special{em:lineto}}
\put(19.63,0.00){\special{em:lineto}}
\put(17.55,0.00){\special{em:lineto}}
\put(17.55,0.00){\special{em:moveto}}
\put(15.63,1.00){\special{em:lineto}}
\put(17.55,2.00){\special{em:lineto}}
\put(17.55,2.00){{\setbox0=\hbox{$\scriptstyle\bullet$}\kern-.4\wd0\lower.5\ht0\box0}}
\put(17.55,0.00){{\setbox0=\hbox{$\scriptstyle\bullet$}\kern-.4\wd0\lower.5\ht0\box0}}
\put(15.63,1.00){{\setbox0=\hbox{$\scriptstyle\bullet$}\kern-.4\wd0\lower.5\ht0\box0}}
\put(21.95,2.00){{\setbox0=\hbox{$\scriptstyle\bullet$}\kern-.4\wd0\lower.5\ht0\box0}}
\put(19.63,0.00){{\setbox0=\hbox{$\scriptstyle\bullet$}\kern-.4\wd0\lower.5\ht0\box0}}
\put(4.16,3.10){\special{em:moveto}}
\put(4.16,4.40){\special{em:lineto}}
\put(4.16,4.40){{\setbox0=\hbox{$\scriptstyle\bullet$}\kern-.4\wd0\lower.5\ht0\box0}}
\put(10.25,3.10){{\setbox0=\hbox{$\scriptstyle\bullet$}\kern-.4\wd0\lower.5\ht0\box0}}
\put(8.16,3.10){{\setbox0=\hbox{$\scriptstyle\bullet$}\kern-.4\wd0\lower.5\ht0\box0}}
\put(6.16,3.10){{\setbox0=\hbox{$\scriptstyle\bullet$}\kern-.4\wd0\lower.5\ht0\box0}}
\put(4.16,3.10){{\setbox0=\hbox{$\scriptstyle\bullet$}\kern-.4\wd0\lower.5\ht0\box0}}
\put(0.16,3.10){{\setbox0=\hbox{$\scriptstyle\bullet$}\kern-.4\wd0\lower.5\ht0\box0}}
\put(0.16,3.10){\special{em:moveto}}
\put(12.25,3.10){\special{em:lineto}}
\put(12.25,7.10){{\setbox0=\hbox{$\scriptstyle\bullet$}\kern-.4\wd0\lower.5\ht0\box0}}
\put(10.16,7.10){{\setbox0=\hbox{$\scriptstyle\bullet$}\kern-.4\wd0\lower.5\ht0\box0}}
\put(2.00,7.10){{\setbox0=\hbox{$\scriptstyle\bullet$}\kern-.4\wd0\lower.5\ht0\box0}}
\put(8.16,7.10){{\setbox0=\hbox{$\scriptstyle\bullet$}\kern-.4\wd0\lower.5\ht0\box0}}
\put(6.08,7.10){{\setbox0=\hbox{$\scriptstyle\bullet$}\kern-.4\wd0\lower.5\ht0\box0}}
\put(4.08,7.10){{\setbox0=\hbox{$\scriptstyle\bullet$}\kern-.4\wd0\lower.5\ht0\box0}}
\put(0.00,7.10){{\setbox0=\hbox{$\scriptstyle\bullet$}\kern-.4\wd0\lower.5\ht0\box0}}
\put(0.00,7.10){\special{em:moveto}}
\put(14.16,7.10){\special{em:lineto}}
\end{picture}}    }
&
\begin{tabular}{c}
\\
$H_1^8$\\
$H_2^8$\\
$H_3^8$\\
$H_4^8$\\
\\
\end{tabular}&

\begin{tabular}{c}
$0.0213042335\times10^{-6}$\\
$0.1810859845\times10^{-6}$\\
$0.3621719690\times10^{-6}$\\
$5.7947515032\times10^{-6}$\\
\end{tabular}\\

\hline

{\scriptsize
\parbox[c]{4.5cm}{
\setlength{\unitlength}{0.06in}
\begin{picture}(28.26,8.41)
\put(13.43,0.00){{\setbox0=\hbox{$\scriptstyle\bullet$}\kern-.4\wd0\lower.5\ht0\box0}}
\put(28.26,4.00){{\setbox0=\hbox{$\scriptstyle\bullet$}\kern-.4\wd0\lower.5\ht0\box0}}
\put(14.88,7.50){{\setbox0=\hbox{$\scriptstyle\bullet$}\kern-.4\wd0\lower.5\ht0\box0}}
\put(13.05,7.50){{\setbox0=\hbox{$\scriptstyle\bullet$}\kern-.4\wd0\lower.5\ht0\box0}}
\special{em:linewidth 0.014in}
\put(0.00,7.50){\special{em:moveto}}
\put(14.88,7.50){\special{em:lineto}}
\put(11.60,0.00){{\setbox0=\hbox{$\scriptstyle\bullet$}\kern-.4\wd0\lower.5\ht0\box0}}
\put(26.43,4.00){{\setbox0=\hbox{$\scriptstyle\bullet$}\kern-.4\wd0\lower.5\ht0\box0}}
\put(3.73,8.41){{\setbox0=\hbox{$\scriptstyle\bullet$}\kern-.4\wd0\lower.5\ht0\box0}}
\put(3.73,7.50){\special{em:moveto}}
\put(3.73,8.41){\special{em:lineto}}
\put(2.36,0.00){{\setbox0=\hbox{$\scriptstyle\bullet$}\kern-.4\wd0\lower.5\ht0\box0}}
\put(22.76,4.00){{\setbox0=\hbox{$\scriptstyle\bullet$}\kern-.4\wd0\lower.5\ht0\box0}}
\put(11.60,1.10){{\setbox0=\hbox{$\scriptstyle\bullet$}\kern-.4\wd0\lower.5\ht0\box0}}
\put(11.60,0.00){\special{em:moveto}}
\put(11.60,1.10){\special{em:lineto}}
\put(21.00,4.00){{\setbox0=\hbox{$\scriptstyle\bullet$}\kern-.4\wd0\lower.5\ht0\box0}}
\put(17.26,4.00){\special{em:moveto}}
\put(17.26,5.10){\special{em:lineto}}
\put(17.26,5.10){{\setbox0=\hbox{$\scriptstyle\bullet$}\kern-.4\wd0\lower.5\ht0\box0}}
\put(27.00,5.58){{\setbox0=\hbox{4}\lower\ht0\box0}}
\put(24.60,4.00){{\setbox0=\hbox{$\scriptstyle\bullet$}\kern-.4\wd0\lower.5\ht0\box0}}
\put(19.10,4.00){{\setbox0=\hbox{$\scriptstyle\bullet$}\kern-.4\wd0\lower.5\ht0\box0}}
\put(17.26,4.00){{\setbox0=\hbox{$\scriptstyle\bullet$}\kern-.4\wd0\lower.5\ht0\box0}}
\put(15.43,4.00){{\setbox0=\hbox{$\scriptstyle\bullet$}\kern-.4\wd0\lower.5\ht0\box0}}
\put(13.60,4.00){{\setbox0=\hbox{$\scriptstyle\bullet$}\kern-.4\wd0\lower.5\ht0\box0}}
\put(13.60,4.00){\special{em:moveto}}
\put(28.26,4.00){\special{em:lineto}}
\put(4.20,0.00){\special{em:moveto}}
\put(4.20,1.18){\special{em:lineto}}
\put(4.20,1.18){{\setbox0=\hbox{$\scriptstyle\bullet$}\kern-.4\wd0\lower.5\ht0\box0}}
\put(9.76,0.00){{\setbox0=\hbox{$\scriptstyle\bullet$}\kern-.4\wd0\lower.5\ht0\box0}}
\put(7.86,0.00){{\setbox0=\hbox{$\scriptstyle\bullet$}\kern-.4\wd0\lower.5\ht0\box0}}
\put(6.03,0.00){{\setbox0=\hbox{$\scriptstyle\bullet$}\kern-.4\wd0\lower.5\ht0\box0}}
\put(4.20,0.00){{\setbox0=\hbox{$\scriptstyle\bullet$}\kern-.4\wd0\lower.5\ht0\box0}}
\put(0.53,0.00){{\setbox0=\hbox{$\scriptstyle\bullet$}\kern-.4\wd0\lower.5\ht0\box0}}
\put(0.53,0.00){\special{em:moveto}}
\put(13.43,0.00){\special{em:lineto}}
\put(11.21,7.50){{\setbox0=\hbox{$\scriptstyle\bullet$}\kern-.4\wd0\lower.5\ht0\box0}}
\put(9.30,7.50){{\setbox0=\hbox{$\scriptstyle\bullet$}\kern-.4\wd0\lower.5\ht0\box0}}
\put(1.83,7.50){{\setbox0=\hbox{$\scriptstyle\bullet$}\kern-.4\wd0\lower.5\ht0\box0}}
\put(7.46,7.50){{\setbox0=\hbox{$\scriptstyle\bullet$}\kern-.4\wd0\lower.5\ht0\box0}}
\put(5.56,7.50){{\setbox0=\hbox{$\scriptstyle\bullet$}\kern-.4\wd0\lower.5\ht0\box0}}
\put(3.73,7.50){{\setbox0=\hbox{$\scriptstyle\bullet$}\kern-.4\wd0\lower.5\ht0\box0}}
\put(0.00,7.50){{\setbox0=\hbox{$\scriptstyle\bullet$}\kern-.4\wd0\lower.5\ht0\box0}}
\end{picture}}    }
&
\begin{tabular}{c}
\\
$H_1^9$\\
$H_2^9$\\
$H_3^9$\\
\\
\end{tabular}&

\begin{tabular}{c}
$0.0004650871\times10^{-7}$\\
$0.1225504411\times10^{-7}$\\
$0.2451008823\times10^{-7}$\\
\end{tabular}\\
\hline
\multicolumn{3}{c}{}\\
\multicolumn{3}{c}{}\\
\end{tabular}

\end{center}

\setcounter{table}{2}

\subsection{Decompositions of the first type}

\subsubsection*{Notation}

Simplices with the dihedral angles $\pi\frac{k}{q}$
are represented by the diagrams similar to the Coxeter diagrams.
A dihedral angle equal to
$\pi\frac{k}{q}$ is represented by a $(q-2)$-fold edge decomposed into
$k$ parts.

The list of the decompositions of the first type
was obtained by the inductive algorithm.
In the tables, the Coxeter diagram of the fundamental simplex is
followed by the non-trivial decompositions with this fundamental simplex.
We omit Coxeter simplices that are not fundamental for non-trivial
decompositions. The non-trivial decompositions are listed in the order
their were obtained in the inductive algorithm.  The tables contain the
simple decompositions only (see section~\ref{type1} for the definition of
the simple decomposition).  The decompositions of the Coxeter simplices
are listed in Table~3.

The numbers ($N,s$ ; $k,l,i,j$) shown under the diagram of the decomposition
describe the decomposition:
\begin{itemize}
\item[]
$\bullet$
$N$ is a number of the fundamental simplices in the decomposition\\
$\bullet$
$s$ is a number of gluing necessary to obtain the decomposition
(if the decomposition is obtained by the gluing together of the simplices
with the gluing numbers
$s_1$ and $s_2$, then $s=1+\max \{s_1,s_2\}$),\\
$\bullet$
$k$ and $l$ are the numbers of the simplices glued together
to obtain the decomposition\\
$\bullet$
$i$ and $j$ are the numbers of the facets of the simplices
$k$ and $l$ which should be glued together
(the nodes of the diagrams are numbered from the left to the right
by the numbers 0,1,...,$n$).

\end{itemize}

\begin{mmtable}
\label{coxtabl}
\begin{center}
\mmcaption{Simple decompositions of the Coxeter simplices of the first type.}
\vspace{7pt}
\begin{tabular}{|c||c|}

\hline
&\\
\begin{tabular}{c}
\begin{tabular}{c|c|c|c|c}
$F$
& $P$ & $N$ & $s$ & $(k,l,i,j)$ \\
\hline

$H_2^{(4)}$ & $H_3^{(4)}$ & 2 & 1 & (0,0,4,4) \\
\hline
$H_6^{4}$ & $H_9^{4}$ & 2 & 1 & (0,0,4,4) \\
$H_7^{4}$ & $H_9^{4}$ & 2 & 1 & (0,0,1,1) \\
$H_5^{4}$ & $H_6^{4}$ & 2 & 1 & (0,0,4,4) \\
$H_5^{4}$ & $H_7^{4}$ & 2 & 1 & (0,0,0,0) \\
$H_4^{4}$ & $H_7^{4}$ & 3 & 2 & (0,0,3,3) \\
$H_4^{4}$ & $H_8^{4}$ & 4 & 2 & (1,0,0,1) \\
$H_2^{4}$ & $H_4^{4}$ & 2 & 1 & (2,2,0,0) \\
$H_2^{4}$ & $H_5^{4}$ & 3 & 2 & (1,0,0,3) \\
$H_1^{4}$ & $H_3^{4}$ & 2 & 1 & (0,0,3,3) \\
$H_1^{4}$ & $H_5^{4}$ & 5 & 3 & (3,2,0,0) \\
\end{tabular}\\
\begin{tabular}{c}
\\
\\
\end{tabular}
\end{tabular}
&
\begin{tabular}{c|c|c|c|c}
$F$ & $P$ & $N$ & $s$ & $(k,l,i,j)$\\
\hline
$H_{10}^{5}$ & $H_{11}^{5}$ & 2 & 1 & (0,0,0,0) \\
$H_8^{5}$    & $H_{10}^{5}$ & 2 & 1 & (0,0,0,0) \\
$H_7^{5}$    & $H_{12}^{5}$ & 8 & 3 & (5,5,0,0) \\
$H_7^{5}$    & $H_{10}^{5}$ & 3 & 2 & (1,0,0,1) \\
$H_5^{5}$    & $H_8^{5}$    & 3 & 2 & (1,0,0,4) \\
$H_5^{5}$    & $H_7^{5}$    & 2 & 1 & (0,0,0,0) \\
$H_3^{5}$    & $H_8^{5}$    & 6 & 3 & (2,2,5,5) \\
$H_3^{5}$    & $H_7^{5}$    & 4 & 3 & (2,0,0,2) \\
$H_2^{5}$    & $H_8^{5}$    & 10 & 4 & (6,6,0,0) \\
$H_2^{5}$    & $H_4^{5}$    & 2 & 1 & (0,0,0,0) \\
$H_1^{5}$    & $H_5^{5}$    & 10 & 4 & (4,6,4,1) \\
$H_1^{5}$    & $H_3^{5}$    & 5 & 4 & (7,0,5,2) \\
$H_1^{5}$    & $H_2^{5}$    & 3 & 2 & (1,0,0,1) \\
\hline
$H_1^{6}$    & $H_2^{6}$    & 2 & 1 & (0,0,5,5) \\
\hline
$H_2^{7}$    & $H_3^{7}$    & 2 & 1 & (0,0,6,6) \\
\hline
$H_2^{8}$    & $H_3^{8}$    & 2 & 1 & (0,0,7,7) \\
\hline
$H_2^{9}$    & $H_3^{9}$    & 2 & 1 & (0,0,8,8) \\

\end{tabular}
\\
&\\
\hline
\end{tabular}

\end{center}
\end{mmtable}

\pagebreak

\begin{mmtable}
\label{pic_hyp1}
\begin{center}
\mmcaption{Simple decompositions of the non-Coxeter simplices of the
first type.}
\end{center}
\end{mmtable}
\vspace{-35pt}
{
\setlength{\unitlength}{0.22175in}

\setlength{\unitlength}{0.2459in}
\begin{center}
{
Decomposition of the bounded simplices in $\H^4$.
}

\end{center}
\vspace{-15pt}

\input{pic-e/simpl/pic_4og.tex}

\setlength{\unitlength}{0.1852in}
\begin{center}
{
Decomposition of the unbounded simplices in $\H^4$.

}
\end{center}
\vspace{-25pt}

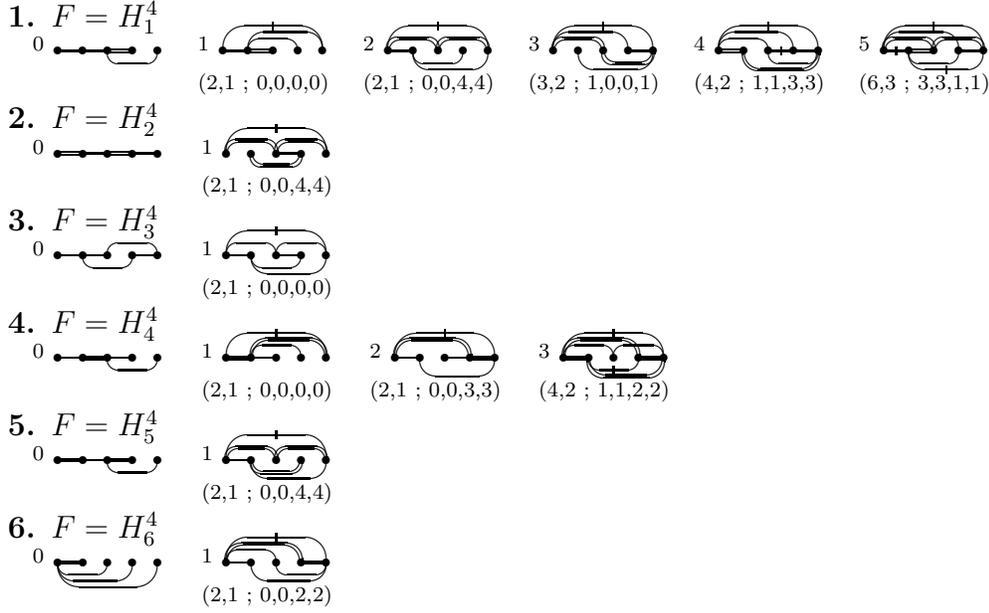
\setlength{\unitlength}{0.1305in}

\begin{picture}(6.35,4.000000)
{\bf \put(-2,1){1. \bf$F=H_{1}^{4}$}}
\multiput(0,0)(1,0){5}{\circle*{0.3}}
\scriptsize\put(-1,0){0}
\put(0,0){\line(1,0){1}}
 \put(1,0){\line(1,0){1}}
 \put(2,-0.05){\line(1,0){1}}
 \put(2,0.05){\line(1,0){1}}
 \put(3.000000,0){\oval(2.000000,1.000000)[b]}
\end{picture}
 \begin{picture}(6.35,4.000000)
\multiput(0,0)(1,0){5}{\circle*{0.3}}
\scriptsize\put(-1,0){1}
{\scriptsize {\put(-1,-1.550000){(2,1 ; 0,0,0,0)}}}
\put(0,0){\line(1,0){1}}
 \put(2.000000,0){\oval(4.000000,2.000000)[t]}
\put(2.000000,0.850000){\line(0,1){0.3}}
\put(1,-0.05){\line(1,0){1}}
 \put(1,0.05){\line(1,0){1}}
 \put(2.000000,0){\oval(2.000000,1.000000)[t]}
\put(2.500000,0){\oval(3.000000,1.500000)[t]}
\end{picture}
 \begin{picture}(6.35,4.000000)
\multiput(0,0)(1,0){5}{\circle*{0.3}}
\scriptsize\put(-1,0){2}
{\scriptsize {\put(-1,-1.550000){(2,1 ; 0,0,4,4)}}}
\put(0,0){\line(1,0){1}}
 \put(1.000000,0){\oval(1.900000,0.900000)[t]}
\put(1.000000,0){\oval(2.100000,1.100000)[t]}
\put(2.000000,0){\oval(4.000000,2.000000)[t]}
\put(2.000000,0.850000){\line(0,1){0.3}}
\put(2.000000,0){\oval(2.000000,1.000000)[b]}
\put(2.500000,0){\oval(3.000000,1.500000)[b]}
\put(3.000000,0){\oval(1.900000,0.900000)[t]}
\put(3.000000,0){\oval(2.100000,1.100000)[t]}
\end{picture}
 \begin{picture}(6.35,4.000000)
\multiput(0,0)(1,0){5}{\circle*{0.3}}
\scriptsize\put(-1,0){3}
{\scriptsize {\put(-1,-1.550000){(3,2 ; 1,0,0,1)}}}
\put(1.000000,0){\oval(1.900000,0.900000)[t]}
\put(1.000000,0){\oval(2.100000,1.100000)[t]}
\put(1.500000,0){\oval(3.000000,1.500000)[t]}
\put(2.000000,0){\oval(4.000000,2.000000)[t]}
\put(2.000000,0.850000){\line(0,1){0.3}}
\put(2.500000,0){\oval(3.000000,1.500000)[b]}
\put(3.000000,0){\oval(1.900000,0.900000)[b]}
\put(3.000000,0){\oval(2.100000,1.100000)[b]}
\put(3,0){\line(1,0){1}}
 \end{picture}
 \begin{picture}(6.35,4.000000)
\multiput(0,0)(1,0){5}{\circle*{0.3}}
\scriptsize\put(-1,0){4}
{\scriptsize {\put(-1,-1.550000){(4,2 ; 1,1,3,3)}}}
\put(0,-0.05){\line(1,0){1}}
 \put(0,0.05){\line(1,0){1}}
 \put(1.000000,0){\oval(2.000000,1.000000)[t]}
\put(1.500000,0){\oval(3.000000,1.500000)[t]}
\put(2.000000,0){\oval(4.000000,2.000000)[t]}
\put(2.000000,0.850000){\line(0,1){0.3}}
\put(2.500000,0){\oval(2.900000,1.400000)[b]}
\put(2.500000,0){\oval(3.100000,1.600000)[b]}
\put(2,0){\line(1,0){1}}
 \put(2.500000,-0.150000){\line(0,1){0.3}}
\put(3.000000,0){\oval(2.000000,1.000000)[b]}
\put(3,0){\line(1,0){1}}
 \end{picture}
 \begin{picture}(6.35,4.000000)
\multiput(0,0)(1,0){5}{\circle*{0.3}}
\scriptsize\put(-1,0){5}
{\scriptsize {\put(-1,-1.550000){(6,3 ; 3,3,1,1)}}}
\put(0,0){\line(1,0){1}}
 \put(0.500000,-0.150000){\line(0,1){0.3}}
\put(1.000000,0){\oval(1.900000,0.900000)[t]}
\put(1.000000,0){\oval(2.100000,1.100000)[t]}
\put(1.500000,0){\oval(3.000000,1.500000)[t]}
\put(2.000000,0){\oval(4.000000,2.000000)[t]}
\put(2.000000,0.850000){\line(0,1){0.3}}
\put(1,-0.05){\line(1,0){1}}
 \put(1,0.05){\line(1,0){1}}
 \put(2.000000,0){\oval(2.000000,1.000000)[b]}
\put(2.500000,0){\oval(3.000000,1.500000)[b]}
\put(2.500000,-0.900000){\line(0,1){0.3}}
\put(3.000000,0){\oval(1.900000,0.900000)[t]}
\put(3.000000,0){\oval(2.100000,1.100000)[t]}
\put(3,0){\line(1,0){1}}
 \end{picture}

\begin{picture}(6.35,4.000000)
{\bf \put(-2,1){2. \bf$F=H_{2}^{4}$}}
\multiput(0,0)(1,0){5}{\circle*{0.3}}
\scriptsize\put(-1,0){0}
\put(0,-0.05){\line(1,0){1}}
 \put(0,0.05){\line(1,0){1}}
 \put(1,0){\line(1,0){1}}
 \put(2,-0.05){\line(1,0){1}}
 \put(2,0.05){\line(1,0){1}}
 \put(3,0){\line(1,0){1}}
 \end{picture}
 \begin{picture}(6.35,4.000000)
\multiput(0,0)(1,0){5}{\circle*{0.3}}
\scriptsize\put(-1,0){1}
{\scriptsize {\put(-1,-1.550000){(2,1 ; 0,0,4,4)}}}
\put(1.000000,0){\oval(1.900000,0.900000)[t]}
\put(1.000000,0){\oval(2.100000,1.100000)[t]}
\put(2.000000,0){\oval(4.000000,2.000000)[t]}
\put(2.000000,0.850000){\line(0,1){0.3}}
\put(2.000000,0){\oval(1.900000,0.900000)[b]}
\put(2.000000,0){\oval(2.100000,1.100000)[b]}
\put(2,0){\line(1,0){1}}
 \put(3.000000,0){\oval(1.900000,0.900000)[t]}
\put(3.000000,0){\oval(2.100000,1.100000)[t]}
\end{picture}

 \begin{picture}(6.35,4.000000)
{\bf \put(-2,1){3. \bf$F=H_{3}^{4}$}}
\multiput(0,0)(1,0){5}{\circle*{0.3}}
\scriptsize\put(-1,0){0}
\put(0,0){\line(1,0){1}}
 \put(1,0){\line(1,0){1}}
 \put(2.000000,0){\oval(2.000000,1.000000)[b]}
\put(3.000000,0){\oval(2.000000,1.000000)[t]}
\put(3,0){\line(1,0){1}}
 \end{picture}
 \begin{picture}(6.35,4.000000)
\multiput(0,0)(1,0){5}{\circle*{0.3}}
\scriptsize\put(-1,0){1}
{\scriptsize {\put(-1,-1.550000){(2,1 ; 0,0,0,0)}}}
\put(0,0){\line(1,0){1}}
 \put(1.000000,0){\oval(2.000000,1.000000)[t]}
\put(2.000000,0){\oval(4.000000,2.000000)[t]}
\put(2.000000,0.850000){\line(0,1){0.3}}
\put(2.000000,0){\oval(2.000000,1.000000)[b]}
\put(2.500000,0){\oval(3.000000,1.500000)[b]}
\put(2,0){\line(1,0){1}}
 \put(3.000000,0){\oval(2.000000,1.000000)[t]}
\end{picture}

 \begin{picture}(6.35,4.000000)
{\bf \put(-2,1){4. \bf$F=H_{4}^{4}$}}
\multiput(0,0)(1,0){5}{\circle*{0.3}}
\scriptsize\put(-1,0){0}
\put(0,0){\line(1,0){1}}
 \put(1,-0.05){\line(1,0){1}}
 \put(1,0.05){\line(1,0){1}}
 \put(2,0){\line(1,0){1}}
 \put(3.000000,0){\oval(2.000000,1.000000)[b]}
\end{picture}
 \begin{picture}(6.35,4.000000)
\multiput(0,0)(1,0){5}{\circle*{0.3}}
\scriptsize\put(-1,0){1}
{\scriptsize {\put(-1,-1.550000){(2,1 ; 0,0,0,0)}}}
\put(0,-0.05){\line(1,0){1}}
 \put(0,0.05){\line(1,0){1}}
 \put(2.000000,0){\oval(4.000000,2.000000)[t]}
\put(2.000000,0.850000){\line(0,1){0.3}}
\put(1,0){\line(1,0){1}}
 \put(2.000000,0){\oval(2.000000,1.000000)[t]}
\put(2.500000,0){\oval(2.900000,1.400000)[t]}
\put(2.500000,0){\oval(3.100000,1.600000)[t]}
\end{picture}
 \begin{picture}(6.35,4.000000)
\multiput(0,0)(1,0){5}{\circle*{0.3}}
\scriptsize\put(-1,0){2}
{\scriptsize {\put(-1,-1.550000){(2,1 ; 0,0,3,3)}}}
\put(0,0){\line(1,0){1}}
 \put(1.500000,0){\oval(2.900000,1.400000)[t]}
\put(1.500000,0){\oval(3.100000,1.600000)[t]}
\put(2.000000,0){\oval(4.000000,2.000000)[t]}
\put(2.000000,0.850000){\line(0,1){0.3}}
\put(2.500000,0){\oval(3.000000,1.500000)[b]}
\put(2,0){\line(1,0){1}}
 \put(3,-0.05){\line(1,0){1}}
 \put(3,0.05){\line(1,0){1}}
 \end{picture}
 \begin{picture}(6.35,4.000000)
\multiput(0,0)(1,0){5}{\circle*{0.3}}
\scriptsize\put(-1,0){3}
{\scriptsize {\put(-1,-1.550000){(4,2 ; 1,1,2,2)}}}
\put(0,-0.05){\line(1,0){1}}
 \put(0,0.05){\line(1,0){1}}
 \put(1.000000,0){\oval(2.000000,1.000000)[t]}
\put(1.500000,0){\oval(2.900000,1.400000)[t]}
\put(1.500000,0){\oval(3.100000,1.600000)[t]}
\put(2.000000,0){\oval(4.000000,2.000000)[t]}
\put(2.000000,0.850000){\line(0,1){0.3}}
\put(2.000000,0){\oval(2.000000,1.000000)[b]}
\put(2.000000,-0.650000){\line(0,1){0.3}}
\put(2.500000,0){\oval(2.900000,1.400000)[b]}
\put(2.500000,0){\oval(3.100000,1.600000)[b]}
\put(3.000000,0){\oval(2.000000,1.000000)[t]}
\put(3,-0.05){\line(1,0){1}}
 \put(3,0.05){\line(1,0){1}}
 \end{picture}

 \begin{picture}(6.35,4.000000)
{\bf \put(-2,1){5. \bf$F=H_{5}^{4}$}}
\multiput(0,0)(1,0){5}{\circle*{0.3}}
\scriptsize\put(-1,0){0}
\put(0,-0.05){\line(1,0){1}}
 \put(0,0.05){\line(1,0){1}}
 \put(1,0){\line(1,0){1}}
 \put(2,-0.05){\line(1,0){1}}
 \put(2,0.05){\line(1,0){1}}
 \put(3.000000,0){\oval(2.000000,1.000000)[b]}
\end{picture}
 \begin{picture}(6.35,4.000000)
\multiput(0,0)(1,0){5}{\circle*{0.3}}
\scriptsize\put(-1,0){1}
{\scriptsize {\put(-1,-1.550000){(2,1 ; 0,0,4,4)}}}
\put(0,0){\line(1,0){1}}
 \put(1.000000,0){\oval(1.900000,0.900000)[t]}
\put(1.000000,0){\oval(2.100000,1.100000)[t]}
\put(2.000000,0){\oval(4.000000,2.000000)[t]}
\put(2.000000,0.850000){\line(0,1){0.3}}
\put(2.000000,0){\oval(1.900000,0.900000)[b]}
\put(2.000000,0){\oval(2.100000,1.100000)[b]}
\put(2.500000,0){\oval(3.000000,1.500000)[b]}
\put(3.000000,0){\oval(1.900000,0.900000)[t]}
\put(3.000000,0){\oval(2.100000,1.100000)[t]}
\end{picture}

 \begin{picture}(6.35,4.000000)
{\bf \put(-2,1){6. \bf$F=H_{6}^{4}$}}
\multiput(0,0)(1,0){5}{\circle*{0.3}}
\scriptsize\put(-1,0){0}
\put(0,-0.05){\line(1,0){1}}
 \put(0,0.05){\line(1,0){1}}
 \put(1.000000,0){\oval(2.000000,1.000000)[b]}
\put(1.500000,0){\oval(3.000000,1.500000)[b]}
\put(2.000000,0){\oval(4.000000,2.000000)[b]}
\end{picture}
 \begin{picture}(6.35,4.000000)
\multiput(0,0)(1,0){5}{\circle*{0.3}}
\scriptsize\put(-1,0){1}
{\scriptsize {\put(-1,-1.550000){(2,1 ; 0,0,2,2)}}}
\put(0,0){\line(1,0){1}}
 \put(1.000000,0){\oval(2.000000,1.000000)[t]}
\put(1.500000,0){\oval(2.900000,1.400000)[t]}
\put(1.500000,0){\oval(3.100000,1.600000)[t]}
\put(2.000000,0){\oval(4.000000,2.000000)[t]}
\put(2.000000,0.850000){\line(0,1){0.3}}
\put(2.500000,0){\oval(3.000000,1.500000)[b]}
\put(3.000000,0){\oval(2.000000,1.000000)[b]}
\put(3,-0.05){\line(1,0){1}}
 \put(3,0.05){\line(1,0){1}}
 \end{picture}
 \\ %

\setlength{\unitlength}{0.215in}
\begin{center}
{
Decomposition of the simplices in $\H^5$.
}
\end{center}
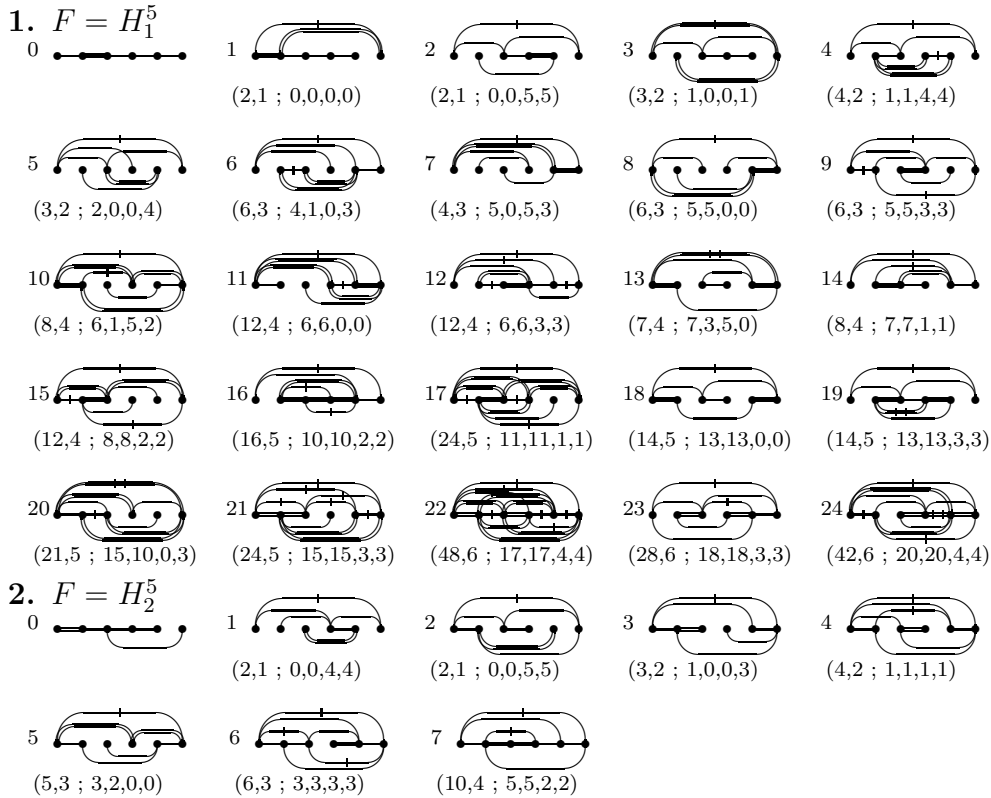
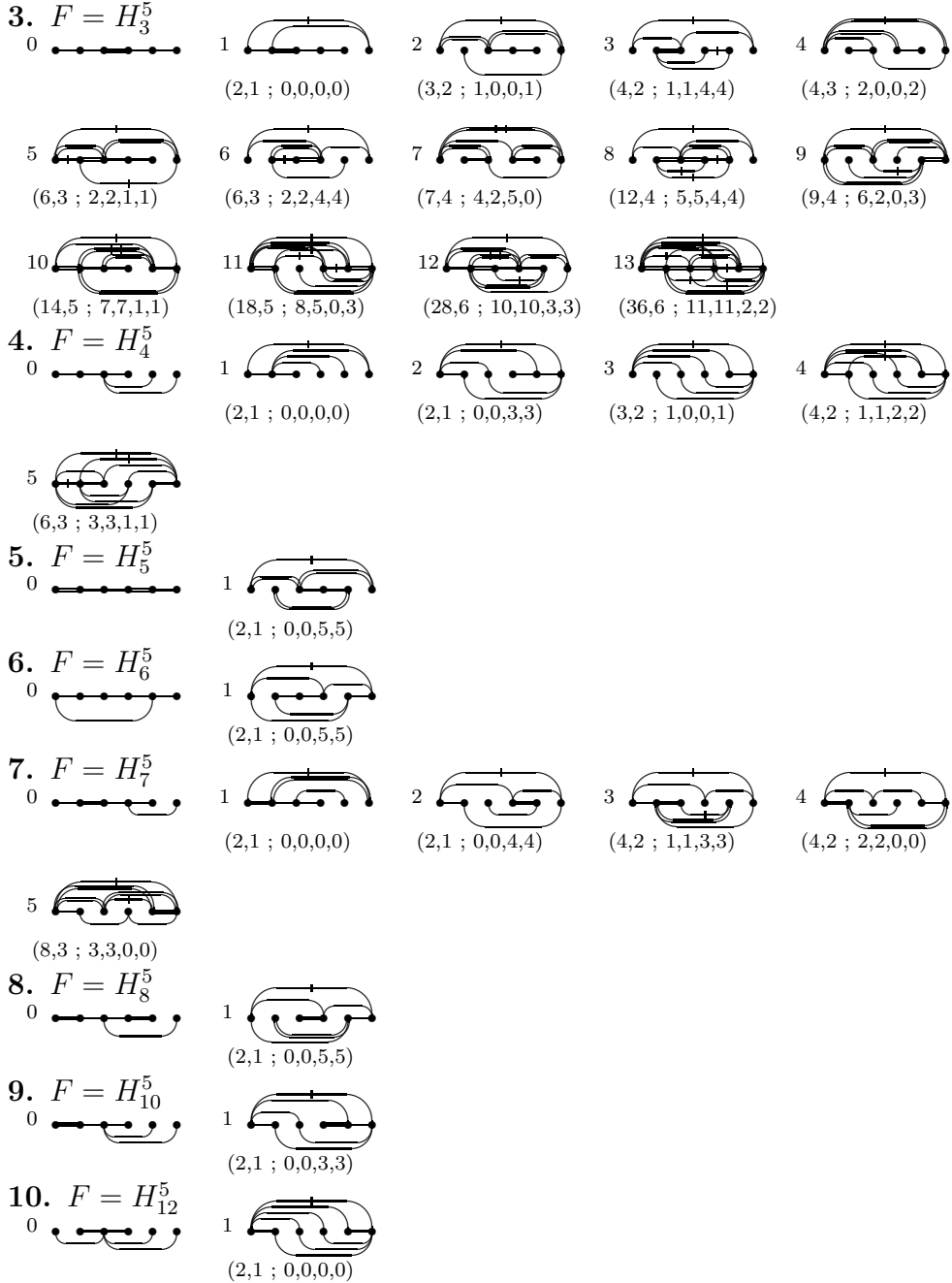
\vspace{-25pt}

\input{pic-e/simpl/pic_5.tex}

\setlength{\unitlength}{0.194in}
\begin{center}
{
Decomposition of the simplices in $\H^6$.
}
\end{center}
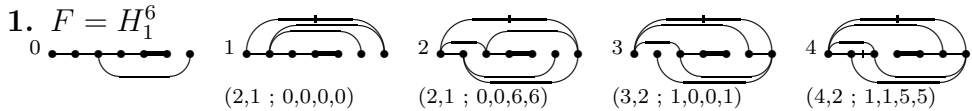
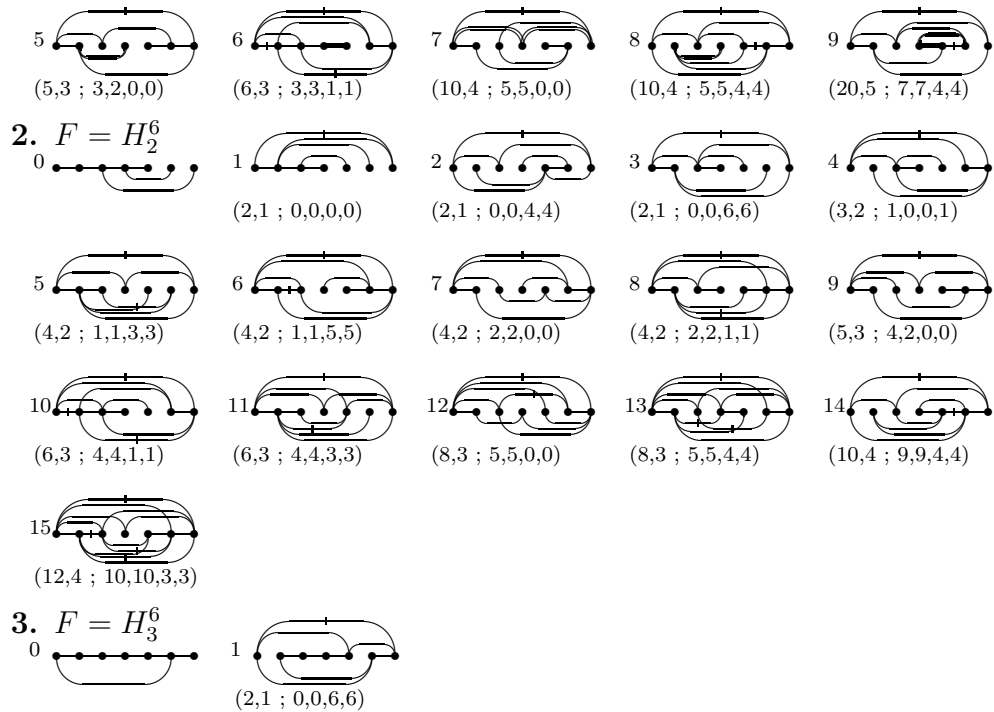
\vspace{-35pt}

\input{pic-e/simpl/pic_6.tex}

\setlength{\unitlength}{0.25in}
\begin{center}
{
Decomposition of the simplices in $\H^7$.
}
\end{center}
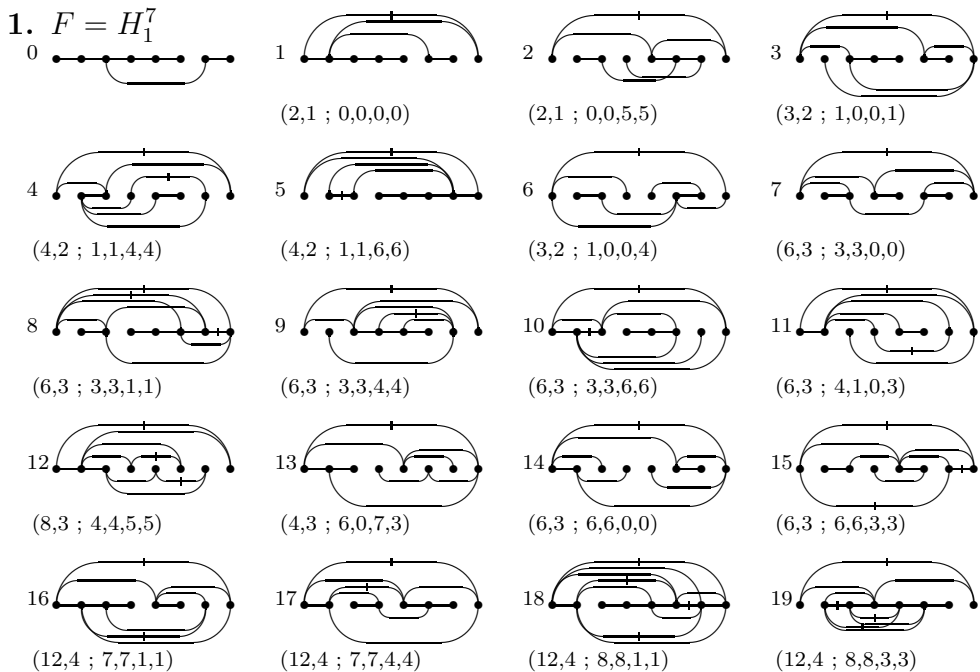
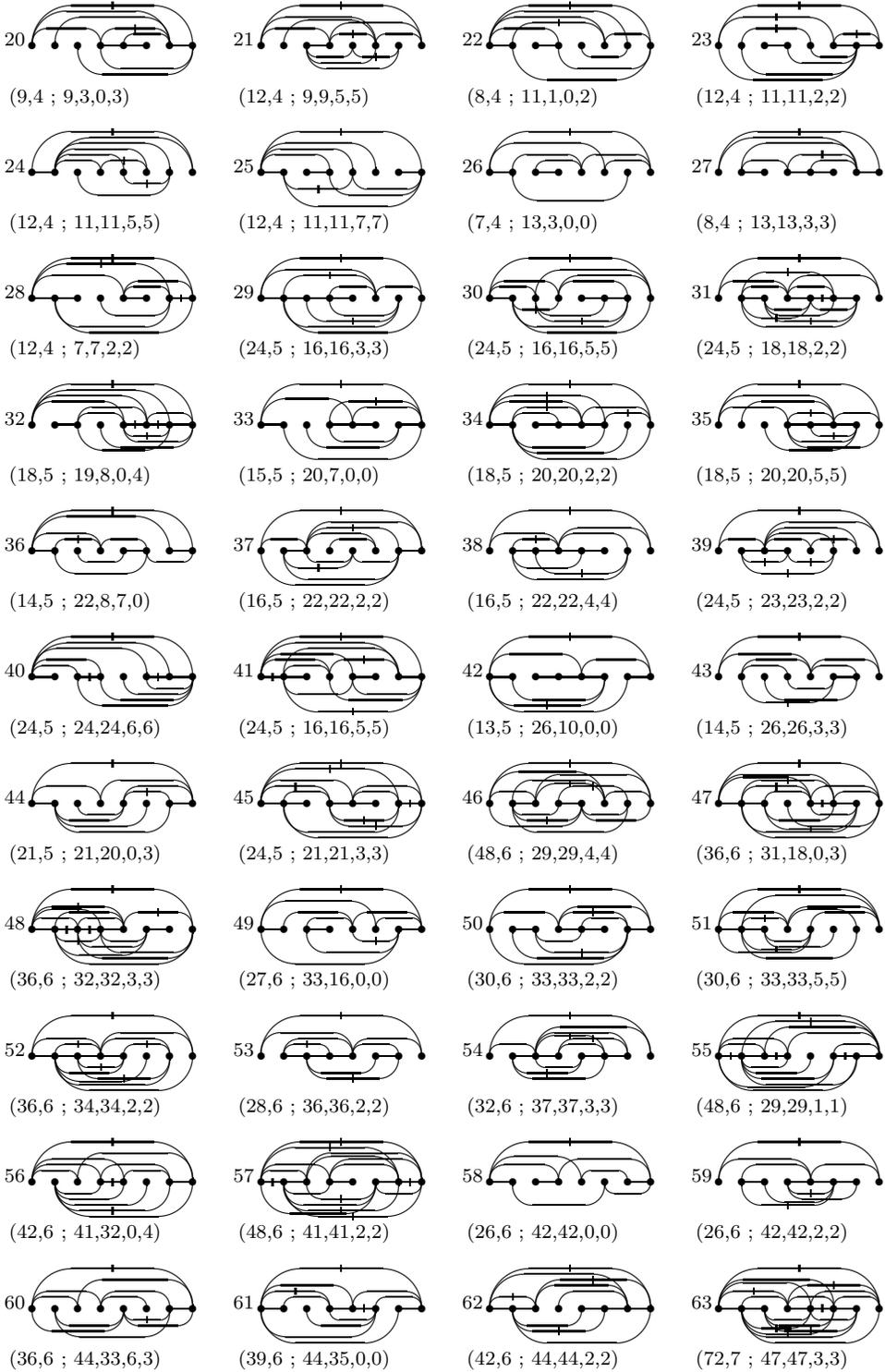
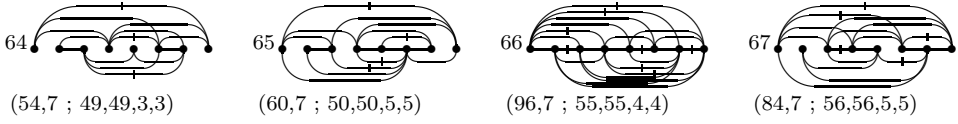
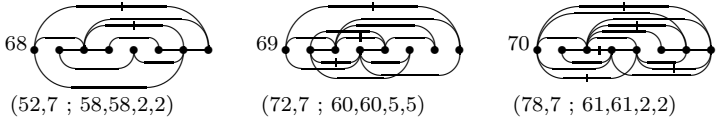
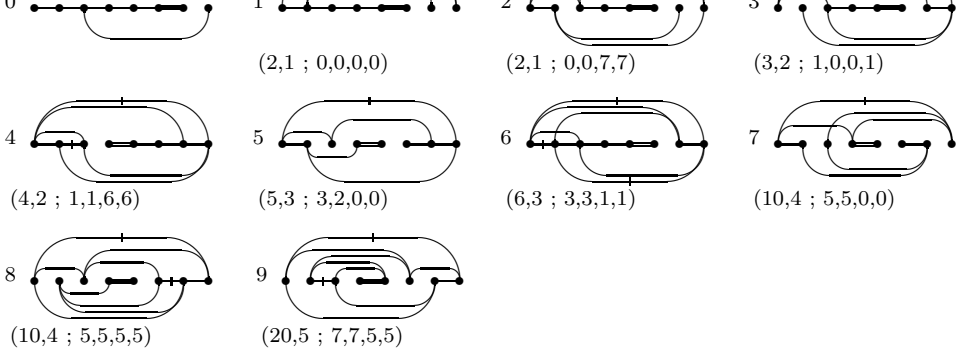
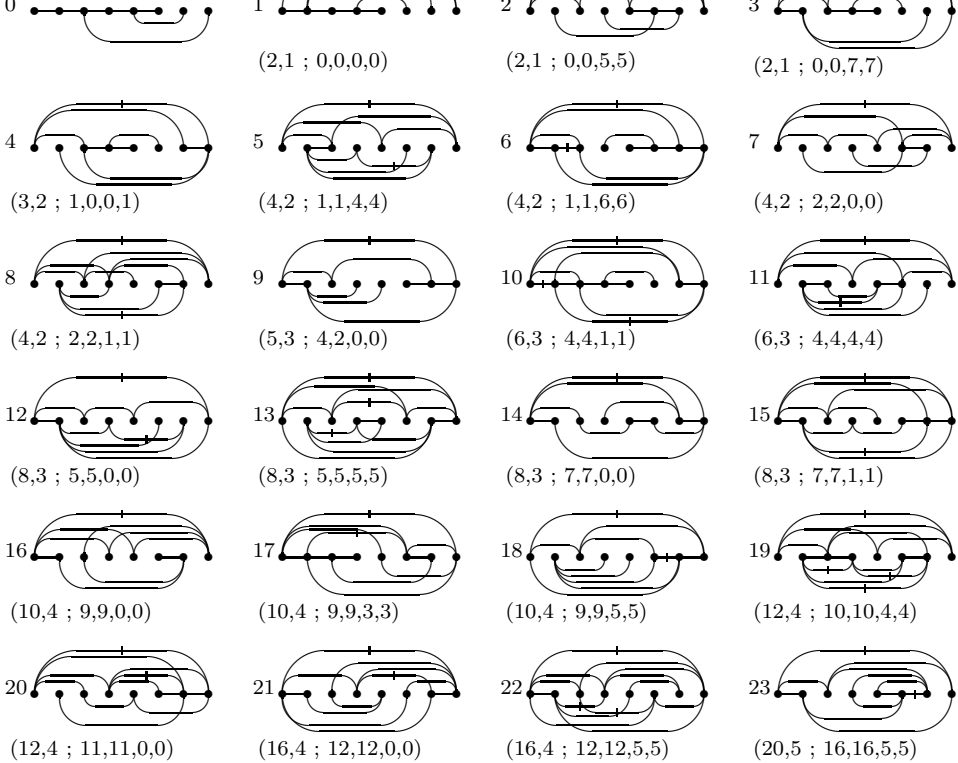
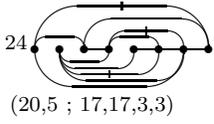
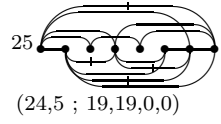
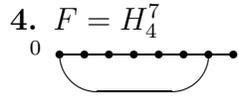
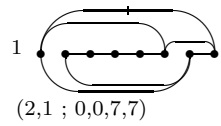
\vspace{-25pt}

\input{pic-e/simpl/pic_7.tex}

\setlength{\unitlength}{0.165in}
\begin{center}
{
Decomposition of the simplices in $\H^8$.

}
\end{center}
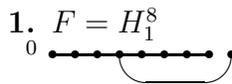
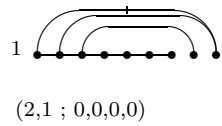
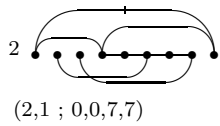
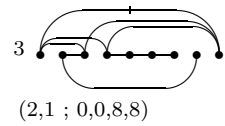
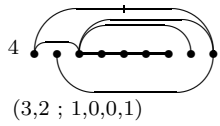
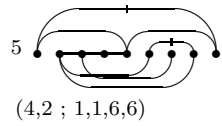
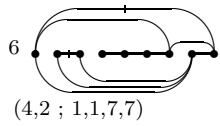
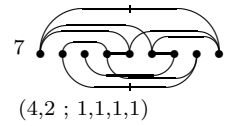
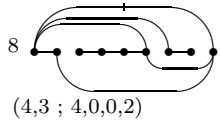
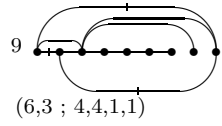
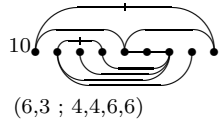
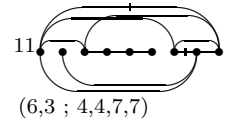
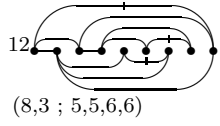
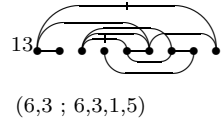
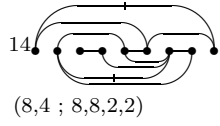
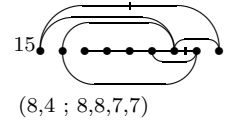
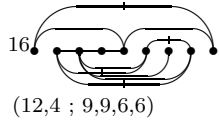
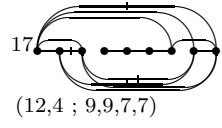
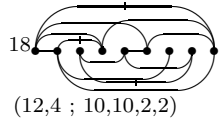
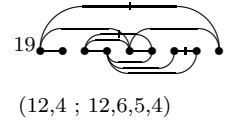
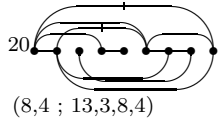
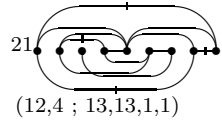
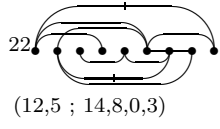
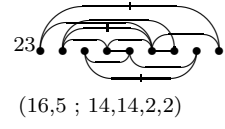
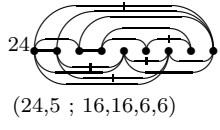
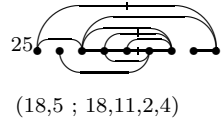
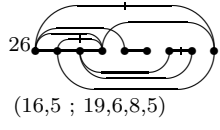
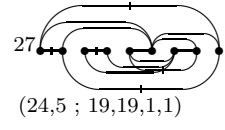
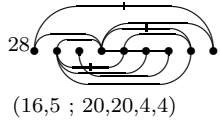
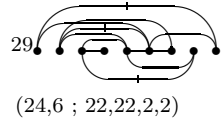
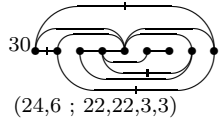
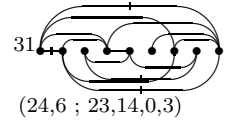
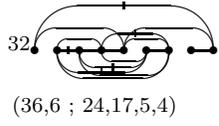
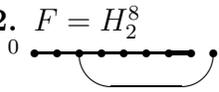
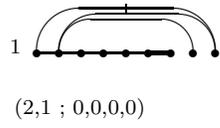
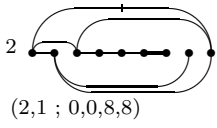
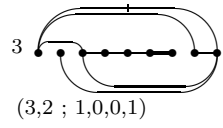
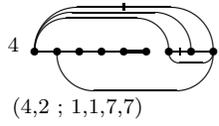
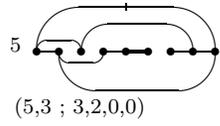
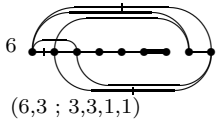
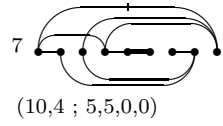
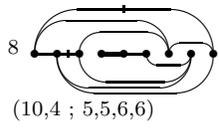
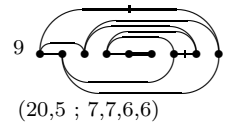
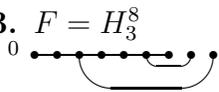
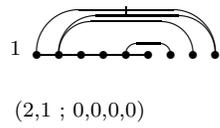
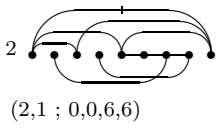
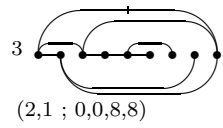
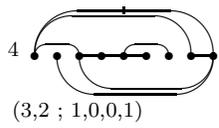
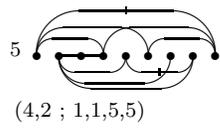
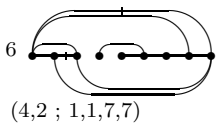
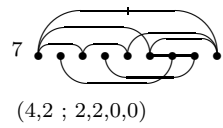
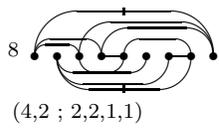
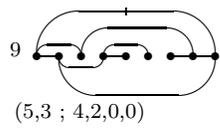
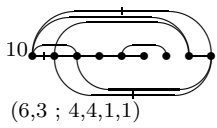
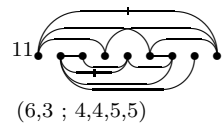
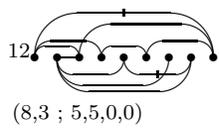
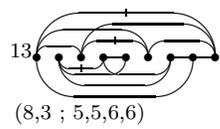
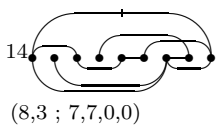
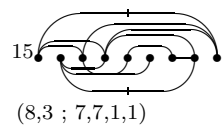
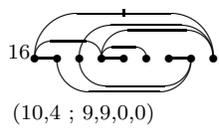
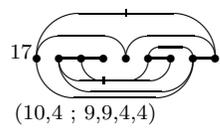
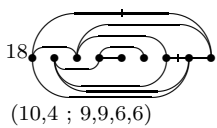
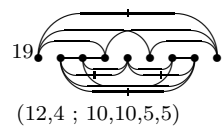
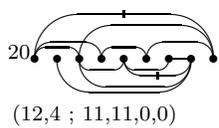
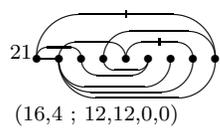
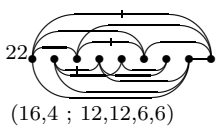
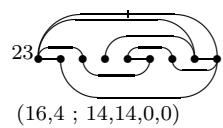
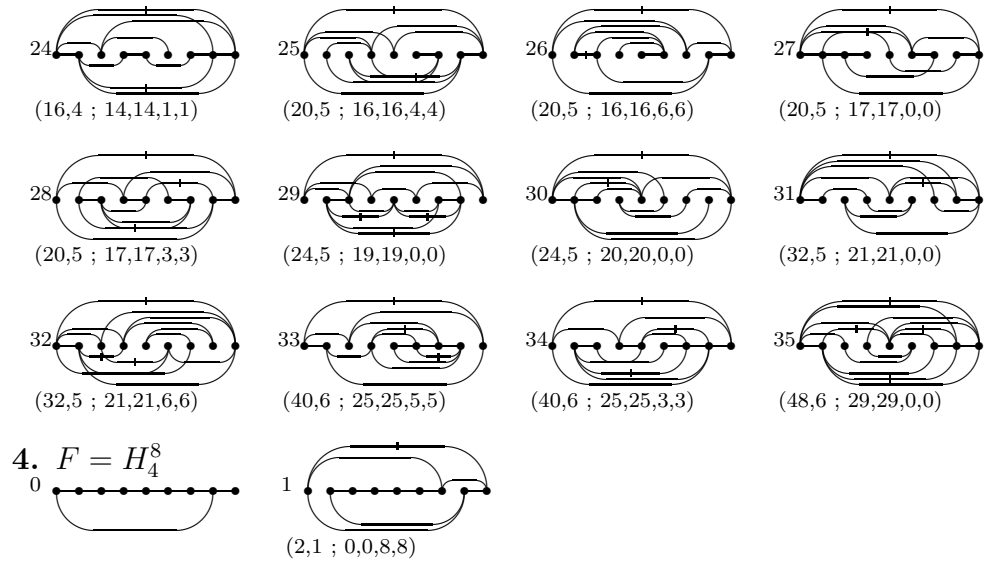
\vspace{-35pt}

\input{pic-e/simpl/pic_8_1.tex}
\setlength{\unitlength}{0.155in}
\vspace{15pt}

\begin{center}

{
Decomposition of the simplices in $\H^9$.
}
\end{center}
\vspace{-39pt}
\input{pic-e/simpl/pic_9_1.tex}

}

\vspace{15pt}

\subsection{Decompositions of the second type}

\begin{mmtable}
\begin{center}
\mmcaption{Decompositions of the second type.}
\vspace{6pt}

\begin{tabular}{|c|c|c|c|}
\hline
$F$ & $P$  &
 $N$ &  Description \\
\hline
$H_3^4$& $H_9^4$&$10$& see section~\ref{type2.4} \\
$H_4^5$& $H_{11}^5$&$20$& see section~\ref{type2.5}\\
$H_1^8$& $H_4^8$&$272$& see section~\ref{type2.8}\\
$H_1^9$& $H_3^9$&$527$& see section~\ref{type2.9}\\
\hline
\end{tabular}

\label{pic_hyp2}
\end{center}
\end{mmtable}



\begin{thebibliography}{13}



\bibitem{Andr}
E.~M.~Andreev. Intersection of plane boundaries of acute-angles polyhedra.
Math. Notes 8 (1971) 761-764.




\bibitem{Deza} A.~Felikson. Coxeter Decompositions of Hyperbolic
Polygons.  European Journal of Combinatorics, (1998) 19,
801--817.

\bibitem{sph_sh}
A.~Felikson.
Coxeter Decompositions of spherical simplices with fundamental
dihedral angles.
Uspehi Mat. Nauk, vol. 57, no.2, 2002 (Russian).



\bibitem{hyp3} A.~Felikson. Coxeter Decompositions of
Hyperbolic Tetrahedra, preprint, Bielefeld, n 98-083.

\bibitem{sph}
A.~Felikson.
Coxeter Decompositions of Spherical Tetrahedra, preprint,
Bielefeld, n 99-053.




\bibitem{pris}
A.~Felikson.
Coxeter Decompositions of Hyperbolic Pyramids
and Triangular Prisms, preprint,
Bielefeld, n 00-006.





\bibitem{Dyn}
E.~B.~Dynkin. Semisimple subalgebras of semisimple Lie algebras,
Amer. Math. Soc. Trans. (2), 6 (1957) 111-244.


\bibitem{Ruth}
N.~W.~Johnson, R.~Kellerhals, J.~G.~Ratcliffe,
S.~T.~Tschantz. The size of a hyperbolic Coxeter simplex,
Transformation  Groups, Vol.4, No 4, 1999, 329--353.

\bibitem{Ruth2}
N.~W.~Johnson, R.~Kellerhals, J.~G.~Ratcliffe,
S.~T.~Tschantz.
Commensurability classes of hyperbolic Coxeter groups.
Linear Algebra and its Applications 345 (2002) 119-147.



\bibitem{Pink} E.~Klimenko, M.~Sakuma. Two-generator discrete
subgroups of ${\it Isom(\mathbf H^2)}$ containing
orientation-reversing elements,
Geom. Dedicata 72 (1998), 247-282.





\bibitem{K}
A.~W.~Knapp. Doubly generated Fuchsian groups, Mich. Math. J.
 1968, v 15, 289--304.



\bibitem{Mat} J.~P.~Matelski.  The
classification of discrete 2--generator subgroups of
$PSL_2(\mathbf R)$, Israel J. Math, 1982, 42,\ 309--317.

\bibitem{M} G.~D.~Mostow. On discontinuous action of monodromy
groups on the complex n-ball, J. of the AMS, 1988, v.1 n 3,
555-586.










\end{thebibliography}
\end{document}